\documentclass{amsart}
\usepackage{amsmath}

\def\-{\mbox{\raisebox{.18cm}{\makebox[0pt]{\hspace{.33cm}.}}$-$}}

\def\s{\\[4pt]}

\def\mr{\mathrm}
\def\tx{\text}

\title{Intuitionistic mathematics and logic}

\author{Joan Rand Moschovakis and Garyfallia Vafeiadou}

\thanks{\\ \mbox{    } This is an English translation, by the authors, of ``Ta enoratika mathimatika kai i logiki tous," published in the Greek volume Stigmes kai Diarkeies, Ed. Anapolitanos, Nefeli 2009. All rights to this translation remain with the authors.
\\ \mbox{    } Research partially supported by a grant from the European Common
Fund and
the Greek Government (PYTHAGORAS II).
\\ \mbox{    }   We want to thank Kostas Mavrommatis, Yiannis Moschovakis and Sifis
Petrakis, who read preliminary versions of the article and gave us their
comments and suggestions.}

\begin{document}

\maketitle

The first seeds of mathematical intuitionism germinated in Europe over a
century ago in the constructive tendencies of Borel, Baire, Lebesgue,
Poincar\'e, Kronecker and others.  The flowering was the work of one man,
Luitzen Egbertus Jan Brouwer, who taught mathematics at the University of
Amsterdam from 1909 until 1951.  By proving powerful theorems on topological
invariants and fixed points of continuous mappings, Brouwer quickly built a
mathematical reputation strong enough to support his revolutionary ideas
about the nature of mathematical activity.  These ideas influenced Hilbert
and G\"odel\footnote{Brouwer's contribution to the dispute over the
foundations of mathematics made a strong impression on all participants.
In particular, his ideas strongly influenced those of his rival Hilbert,
as well as the work of G\"odel.  For a full discussion of the history
see \cite{Hs2003}.} and established intuitionistic logic and
mathematics as subjects worthy of independent study.

\vskip 0.1cm

Our aim is to describe the development of Brouwer's intuitionism, from
his rejection of the classical law of excluded middle to his controversial
theory of the continuum, with fundamental consequences for logic and
mathematics.  We borrow Kleene's formal axiomatic systems (incorporating
earlier attempts by Kolmogorov, Glivenko, Heyting and Peano) for intuitionistic
logic and arithmetic as subtheories of the corresponding classical theories,
and sketch his use of g\"odel numbers of recursive functions to {\em realize}
sentences of intuitionistic arithmetic including a form of Church's Thesis.
Finally, we present Kleene and Vesley's axiomatic treatment of Brouwer's
continuum, with the function-realizability interpretation which establishes
its consistency.

\vskip 0.1cm

\section{Brouwer's Early Philosophy}

In 1907 L. E. J. Brouwer published (in Dutch) his doctoral dissertation, whose
title can be translated ``On the Foundations of Mathematics.''\footnote{An
English translation, from which the quotations in this section are taken,
appears in \cite{Br1907}.}  This remarkable manifesto, with its heterodox
views on mathematics, logic and language, critically examined and found fault
with every major mathematical philosophy of the time.  While Brouwer was aware
of the work of the French intuitionist Poincar\'e and of Borel's constructive
approach to the theory of sets, it was not in his nature to be a follower.
In his 1912 inaugural address at the University of Amsterdam he referred to
his own philosophy as {\em neo-intuitionism},\footnote{``Intuitionism and
Formalism,'' an English translation of this address, appeared in the Bulletin
of the American Mathematical Society in the same year, and is included in
\cite{BP1964}.} but the vigor and creativity Brouwer brought to the subject
over almost half a century of work have linked his name, more than any other,
with intuitionistic philosophy and mathematics.  Many of the basic principles
of his intuitionism were already clear in his dissertation.

\vskip 0.1cm

In direct opposition to Russell and Whitehead's logicism, Brouwer
asserted in 1907 that mathematics cannot be considered a part of logic.
``Strictly speaking the construction of intuitive mathematics in itself
is an {\em action} and not a {\em science}; it only becomes a science,
i.e. a totality of causal sequences, repeatable in time, in a mathematics
of the second order [metamathematics], which consists of the {\em mathematical
consideration of mathematics} or {\em of the language of mathematics}...
But there, as in the case of theoretical logic, we are concerned with an
{\em application of mathematics}, that is, with an {\em experimental science}''
(\cite{Br1907} p. 61).\footnote{Within quotations all words in [~] are ours,
but the italics are Brouwer's own.}

\vskip 0.1cm

The discovery of non-Euclidean geometries showed, according to Brouwer,
that Kant was only partly right in asserting that the intuitions of space
and time are logically prior to (and independent of) experience.  ``...we
can call a priori only that one thing which is common to all mathematics
and is \ldots sufficient to build up all mathematics, namely the intuition
of the many-oneness, the basic intuition of mathematics.   And since in this
intuition we become conscious of time as change per se, we can state:
{\em the only a priori element in science is time}'' (\cite{Br1907} p. 61).

\vskip 0.1cm

Hilbert's formalist program was doomed to failure because ``language ... is
a means ... for the communication of mathematics but ... has nothing to do
with mathematics'' and is not essential for it.  Moreover, the
``... {\em existence} of a mathematical system satisfying a set of axioms
can never be proved from the consistency of the logical system based on those
axioms,'' but only by construction.  ``A fortiori it is not certain that any
mathematical problem can either be solved or proved to be unsolvable''
(\cite{Br1907} p. 79).

\vskip 0.1cm

According to Brouwer, the paradoxes in set theory arise from the
consideration of sets which are too large and abstract to be built up
mathematically.  Even Cantor's second number class (of the denumerably
infinite ordinals) cannot exist, although the concept is consistent.
Zermelo's proof of the wellordering principle from the axiom of choice
is misguided.  The continuum cannot be well-ordered, ``firstly because
the greater part of the elements of the continuum must be considered as
unknown, and ... secondly because every well-ordered set is denumerable''
(\cite{Br1907} pp. 84-85).

\vskip 0.1cm

At the end of the dissertation, as the second of twenty-one ``STATEMENTS
(to be defended together with the thesis),'' Brouwer asserts:  ``It is not
only impossible to prove the admissibility of complete induction, but it
ought neither to be considered as a special axiom nor as a special intuitive
truth.  Complete induction is an act of mathematical construction, justified
simply by the basic intuition of mathematics'' (\cite{Br1907} p. 98).
This statement effectively dismisses the work of Peano, but admits a potential
infinity of natural numbers with a method for showing that arbitrarily
complicated properties (even those involving quantification over all natural
numbers) hold for each.  In this context Brouwer's earlier remark that
``{\em all} or {\em every} \ldots tacitly involves the restriction:
{\em insofar as belonging to a mathematical structure which is supposed to be
constructed beforehand}'' (\cite{Br1907} p. 76) suggests that the {\em structure}
of the natural numbers can be understood as a completed construction even though
the {\em collection} of all natural numbers cannot be surveyed at a glance.

\vskip 0.1cm

\section{Intuitionistic Logic}

One year after his dissertation, in ``The unreliability of the logical
principles,'' Brouwer argued against the use of classical logic in mathematics
and science.  He agreed with the principles of {\em syllogism} (if all A's are B,
and all B's are C, then all A's are C) and {\em contradiction} (nothing is both
A and not A), but not with the law of excluded middle (everything is A or not A)
when it is applied to infinite systems.\footnote{\cite{Br1908} pp. 109-110.
In fact, Brouwer asserted, ``\ldots the question of the validity of the principium
tertii exclusi is equivalent to the question {\em whether unsolvable mathematical
problems can exist}.  There is not a shred of a proof for the conviction, which
has sometimes been put forward, that there exist no unsolvable mathematical
problems.'' We shall return to this question in Sections 3 and 4.} In effect, Brouwer
distinguished between the intuitionistically unacceptable $A \vee \neg A$
and the intuitionistically correct $\neg \neg (A \vee \neg A)$ , making
full use of the expressive power of the logical language to separate constructions
which establish facts from constructions which establish consistency.\footnote{D.
Hesseling (\cite{Hs2003} p.280) credits A. Kolmogorov \cite{Ko1932} with the
observation that for Brouwer a statement of the form $\neg B$ was {\em positive
existential}: ``there exists a chain of logical inferences, starting with the
assumption that [$B$] is correct and concluding with a contradiction'' .  By the
same reasoning, $\neg \neg B$ asserts the consistency of $B$.}

\vskip 0.1cm

Considering Brouwer's attitude toward formal logic, it is hardly surprising
that he did not attempt to axiomatize intuitionistic reasoning.  Nevertheless,
he recognized the usefulness of formulating general principles which could
be relied on for mathematical constructions.  It is this which legitimizes
formal systems for intuitionistic logic and mathematics, as long as each
axiom and rule can be justified from an intuitionistic standpoint.
As Kleene remarks on page 5 of \cite{FIM}, Brouwer's objection was only to
formal reasoning without a corresponding (mathematical) meaning.

\vskip 0.1cm

In 1925 Andrei Kolmogorov \cite{Ko1925} proposed axioms for (minimal)
intuitionistic reasoning with implication and negation only; his article,
in Russian, attracted little or no attention in western Europe.  Formal systems
for intuitionistic propositional logic were published (in French) by Valerii
Glivenko in 1928 (\cite{Gl1928}) and 1929 (\cite{Gl1929}).  The first was
incomplete.  The second included two additional axioms suggested by Brouwer's
student Arend Heyting, who presented his own detailed axiomatizations of
intuitionistic propositional and predicate logic and a part of intuitionistic
mathematics in three classic papers \cite{He1930A}, \cite{He1930B}, \cite{He1930C}
the following year.\footnote{Glivenko's 1929 note and the first part of Heyting's
axiomatization were recently translated from the original French and German into
English for the collection \cite{Ma1998}.}  Heyting began with the caveat
\begin{quote}
``Intuitionistic mathematics is a mental process, and every language, the
formalistic one included, is an aid to communication only.  It is in principle
impossible to construct a system of formulas equivalent to intuitionistic
mathematics, since the possibilities of thinking cannot be reduced to a finite
number of rules constructed in advance.''\footnote{An authoritative and readable
summary of Heyting's contributions to the subject is A. S. Troelstra's
\cite{Tr1968}, from which this smooth English translation of the introduction to
Heyting \cite{He1930A} is borrowed.}
\end{quote}

Nevertheless, the effort to explain in metamathematical terms the difference
between intuitionistic and classical logic made sense, because Brouwer's
theory of the continuum contradicted classical logic.

\vskip 0.1cm

In order to be acceptable intuitionistically, a general logical principle
must be uniformly interpretable by means of constructions (proofs or
computations).  Kolmogorov's problem interpretation, and Heyting's association
of assertions with their (constructively acceptable) proofs, are special cases
of what has become known as the Brouwer-Heyting-Kolmogorov explication of the
intuitionistic connectives and quantifiers.

\vskip 0.1cm

Like Tarski's truth definition,
the B-H-K interpretation is heuristic rather than mathematically precise,
and is based on the assumption that true prime statements justify themselves.
An implication can only be justified by a construction which uniformly converts
any given justification of the hypothesis to a justification of the conclusion.
A disjunction is justified by a construction which chooses a particular one of
the disjuncts and provides its justification.  The negation of a statement is
justified by a construction which would convert any justification of the statement
into a proof of a  known contradiction.  The intuitionistic objections to the
classical laws of excluded middle and double negation can be explained on
the basis of this interpretation.

\vskip 0.1cm

Beginning around 1940 the American logician S. C. Kleene, who was sympathetic
to intuitionistic ideas, devoted considerable effort to clarifying their
precise relation to classical logic and mathematics.\footnote{ Kleene was
inspired by a remark of Hilbert and Bernays to consider ``$\exists x A(x)$,''
for example, as ``an incomplete communication, which is completed by giving
an $x$ such that $A(x)$ together with the further information required to
complete the communication `$A(x)$' for that $x$.''  This interpretation
led to recursive realizability; see \S 5.}  The Hilbert-style formalisms for
intuitionistic logic and arithmetic we present here are from Kleene \cite{IM}
(p. 82), with minor changes in symbols and numbering.  They represent an
independent selection among many axioms proposed earlier by Kolmogorov,
Glivenko, Heyting and Peano.

\subsection{The intuitionistic propositional calculus Pp}

The first step in the metamathematical study of any part of logic or
mathematics is to specify a {\em formal language}.  For propositional or
sentential logic, the standard language has denumerably many distinct
proposition letters $\mathrm{P}_0, \mathrm{P}_1, \mathrm{P}_2, \ldots$ and symbols
 $\&$, $\vee$, $\rightarrow$, $\neg$ for the propositional connectives ``and,'' ``or,''
``if \ldots then,'' and ``not'' respectively, with left and right parentheses
(, ) (sometimes written ``[, ]'' for ease of reading).  Classical logic
actually needs only two connectives (since classical $\vee$ and $\rightarrow$
can be defined in terms of $\&$ and $\neg$), but the four intuitionistic
connectives are independent.  The classical language is thus properly contained
in the intuitionistic, which is more expressive.

\vskip 0.1cm

The most important tool of metamathematics is generalized induction, a method
Brouwer endorsed.  The class of ({\it well-formed}\,) {\it formulas} of the
language of {\bf Pp} is defined inductively by
\begin{enumerate}

\item[(i)]{Each proposition letter is a ({\it prime}\,) {\it formula}.}

\item[(ii)]{If $\mathrm{A}$, $\mathrm{B}$ are {\it formulas} so are $\mathrm{(A ~\&~ B)}$,
$\mathrm{(A \vee B)}$, $\mathrm{(A \rightarrow B)}$ and $\mathrm{(\neg A)}$.}

\item[(iii)]{Nothing is a {\it formula} except as required by (i) and (ii).}
\end{enumerate}
As in classical logic, $\mathrm{(A \leftrightarrow B)}$ abbreviates
$\mathrm{((A \rightarrow B)~\&~(B \rightarrow A))}$.  Inessential parentheses are
omitted on the convention that $\neg$ binds closer than $\&$, $\vee$ which bind closer
than $\rightarrow$.

\vskip 0.1cm

The building blocks for Kleene's version of intuitionistic propositional logic
{\bf Pp} are finitely many {\em axiom schemas}, each summarizing a potentially
infinite collection of intuitionistically correct formulas, and one {\em rule
of inference} expressing an intuitionistically acceptable principle of reasoning
from hypotheses to a conclusion.  The axioms are all formulas of the following
forms:\footnote{Glivenko's original axiom system consisted of X4 - X9 and variants
of X2, X3.  In \cite{Gl1929} he added X1 and X10, which he attributed to Heyting.
Heyting \cite{He1930A} had X1, X6, X9, X10 and variants of X2 and X8 as axioms but
otherwise $\&$ and $\vee$ were treated differently. Heyting showed that his axioms are
independent and do not prove $\mathrm{\neg \neg A \rightarrow A}$.  Kleene's version
of X2 was designed to simplify the proof of the Deduction Theorem.}
\begin{enumerate}
\item[X1.]{$\mathrm{A \rightarrow (B \rightarrow A)}$.}

\item[X2.]{$\mathrm{(A \rightarrow B) \rightarrow ((A \rightarrow (B \rightarrow C))
\rightarrow (A \rightarrow C))}$.}

\item[X3.]{$\mathrm{A \rightarrow (B \rightarrow A ~\&~ B)}$.}

\item[X4.]{$\mathrm{A ~\&~ B \rightarrow A}$.}

\item[X5.]{$\mathrm{A ~\&~ B \rightarrow B}$.}

\item[X6.]{$\mathrm{A \rightarrow A \vee B}$.}

\item[X7.]{$\mathrm{B \rightarrow A \vee B}$.}

\item[X8.]{$\mathrm{(A \rightarrow C) \rightarrow ((B \rightarrow C)
\rightarrow (A \vee B \rightarrow C))}$.}

\item[X9.]{$\mathrm{(A \rightarrow B) \rightarrow ((A \rightarrow \neg B)
\rightarrow \neg A)}$.}

\item[X10.]{$\mathrm{\neg A \rightarrow (A \rightarrow B)}$.}
\end{enumerate}
The rule of inference of {\bf Pp} is
\begin{enumerate}
\item[R1]{({\it Modus Ponens}).  From $\mathrm{A}$ and $\mathrm{A \rightarrow B}$,
conclude $B$.}
\end{enumerate}

\vskip 0.1cm

A formal {\em proof~} in {\bf Pp} is any finite sequence $\mathrm{E}_1,\ldots,\mathrm{E}_k$
of formulas, each of which is an axiom or an immediate consequence, by the
rule of inference, of two preceding formulas of the sequence.\footnote{This
and similar descriptions abbreviate the obvious inductive definitions.}
Any proof is said to {\it prove} its last formula, which is therefore a
{\it theorem} of {\bf Pp}.  We write $\vdash_{\bf Pp} \mathrm{E}$  to denote that
$\mathrm{E}$ is a theorem of {\bf Pp}.

\vskip 0.1cm

{\it Example}.  Here is a complete formal proof (really a proof schema, as
$\mathrm{A}$ may be any formula) in {\bf Pp} of $\mathrm{\neg (A ~\&~ \neg A)}$,
indicating the reasons for each step.
\begin{enumerate}
\item{$\mathrm{(A ~\&~ \neg A) \rightarrow A}$. [axiom by X4]}

\item{$\mathrm{(A ~\&~ \neg A) \rightarrow \neg A}$. [axiom by X5]}

\item{$\mathrm{((A ~\&~ \neg A) \rightarrow A) \rightarrow
(((A ~\&~ \neg A) \rightarrow \neg A) \rightarrow \neg(A ~\&~ \neg A))}$. [axiom by X9]}

\item{$\mathrm{((A ~\&~ \neg A) \rightarrow \neg A) \rightarrow \neg (A ~\&~ \neg A)}$.
[by R1 from (1), (3)]}

\item{$\mathrm{\neg (A ~\&~ \neg A)}$. [by R1 from (2), (4)]}
\end{enumerate}

\vskip 0.1cm

If $\Gamma$ is any collection of formulas and $~\mathrm{E}_1,\ldots,\mathrm{E}_k~$ any
finite sequence of formulas each of which is a member of $\Gamma$, an axiom, or an
immediate consequence by R1 of two preceding formulas, then
$~\mathrm{E}_1,\ldots,\mathrm{E}_k~$ is a {\it derivation in} {\bf Pp} {\it of} its
last formula $\mathrm{E}_k$ {\it from} the assumptions $\Gamma$.  We write
$\Gamma \vdash_{\bf Pp} \mathrm{E}$ to denote that such a derivation exists with
$\mathrm{E}_k = \mathrm{E}$.  The following theorem is proved by induction
over the definition of a derivation; its converse follows from R1.

\vskip 0.1cm

{\it The Deduction Theorem}.  If $\Gamma$ is any collection of formulas
and $\mathrm{A, B}$ are any formulas such that
$\Gamma \cup \{\mathrm{A}\} \vdash_{\bf Pp} \mathrm{B}$,
then also $\Gamma \vdash_{\bf Pp} \mathrm{(A \rightarrow B)}$.

\vskip 0.1cm

The axiomatization is designed so that classical propositional logic
{\bf Pp$^c$} results from {\bf Pp} by strengthening the axiom schema X10 to
\begin{enumerate}
\item[X10$^c$.]{$\mathrm{\neg \neg A \rightarrow A}$.}
\end{enumerate}
The definitions of {\it proof} and {\it derivation} for {\bf Pp$^c$} are like
those for {\bf Pp} but with X1-X10$^c$ instead of X1-X10.  To show that {\bf Pp}
is a subtheory of {\bf Pp$^c$} it suffices to prove $\vdash_{\bf Pp^c}$ X10,
a relatively simple exercise.

\vskip 0.1cm

In 1929 Glivenko \cite{Gl1929} proved that if $\mathrm{A}$ is any formula such that
$~\vdash_{\bf Pp^c} \mathrm{A}~$, then $~\vdash_{\bf Pp} \mathrm{\neg \neg A}\,$.  This
simple form holds {\em only} for propositional logic, and is known as ``Glivenko's
Theorem.''

\vskip 0.1cm

Around 1933 Kurt G\"odel \cite{Go1933} and Gerhard Gentzen (published
posthumously in \cite{Ge1965}) observed independently that {\bf Pp$^c$}
can be faithfully translated into {\bf Pp}.\footnote{In 1925 Kolmogorov
\cite{Ko1925} published a different negative translation for a fragment
with $\rightarrow$ and $\neg$ only.}  Briefly, each proposition letter
is replaced by its double negation, and then $\vee$ is replaced inductively
by its classical definition in terms of $\neg$ and $\&$.  If $\Gamma^g, \mathrm{A}^g$
are the translations of $\Gamma, \mathrm{A}$ respectively then
\begin{enumerate}
\item[(i)]{$\vdash_{\bf Pp^c} (\mathrm{A}^g \leftrightarrow \mathrm{A})$, and}
\item[(ii)]{$\Gamma^g \vdash_{\bf Pp} \mathrm{A}^g$ if and only if
$\Gamma \vdash_{\bf Pp^c} \mathrm{A}$.}
\end{enumerate}

In 1934-35 Gentzen \cite{Ge1935} proved a normal form theorem for an
intuitionistic sequent calculus, giving an effective algorithm for deciding
whether an arbitrary formula $A$ is or is not provable in {\bf Pp}.  Since
intuitionistic propositional logic has no finite truth-table interpretation,
the decision algorithm for {\bf Pp} is more complicated than for {\bf Pp$^c$}.

\subsection{The intuitionistic first-order predicate calculus Pd}

The pure first-order language has individual variables
$\mathrm{a}_1,\mathrm{a}_2,\mathrm{a}_3,\ldots$,
and countably infinitely many distinct predicate letters
$\mathrm{P}_1(\ldots)$,$\mathrm{P}_2(\ldots)$,$\mathrm{P}_3(\ldots),\ldots$ of arity $n$
for each $n=0,1,2,\ldots$, including the 0-ary proposition letters.  There are
two new logical symbols $\forall$ (``for all'') and $\exists$ (``there exists'').

\vskip 0.1cm

The {\it terms} of the language of {\bf Pd} are the individual variables.
The {\em formulas} are defined by
\begin{enumerate}
\item[(i)]{If $\mathrm{P(\ldots)}$ is an $n$-ary predicate letter and
$\mathrm{t}_1,\ldots,\mathrm{t}_n$
are terms then $\mathrm{P(t}_1,\ldots,\mathrm{t}_n)$ is a ({\em prime}) {\em formula}.}

\item[(ii)]{If $\mathrm{A}$, $\mathrm{B}$ are {\em formulas} so are $\mathrm{(A ~\&~ B)}$,
$\mathrm{(A \vee B)}$, $\mathrm{(A \rightarrow B)}$ and $\mathrm{(\neg A)}$.}

\item[(iii)]{If $\mathrm{A}$ is a {\em formula} and $\mathrm{x}$ an individual variable,
then $\mathrm{(\forall x A)}$ and $\mathrm{(\exists x A)}$ are {\em formulas}.}

\item[(iv)]{Nothing else is a {\em formula}.}
\end{enumerate}

We use $\mathrm{x,y,z,w,x}_1,\mathrm{y}_1,\ldots$ and
$\mathrm{A,B,C,\ldots,A(x),A(x,y),\ldots}$ as metavariables for variables and formulas,
respectively.  Anticipating applications (e.g. to arithmetic),
$\mathrm{s,t,s}_1,\mathrm{t}_1,\ldots$ vary over terms.  In omitting parentheses,
$\mathrm{\forall x}$ and $\mathrm{\exists x}$ are treated like $\neg$.
The {\it scope} of a quantifier, and {\em free} and {\em bound} occurrences
of a variable in a formula, are defined as usual.  A formula in which every
variable is bound is a {\em sentence} or {\em closed formula}.

\vskip 0.1cm

If $\mathrm{x}$ is a variable, $\mathrm{t}$ a term, and $\mathrm{A(x)}$ a formula which
may or may not contain $\mathrm{x}$ free, then $\mathrm{A(t)}$ denotes the result of
substituting an occurrence of $\mathrm{t}$ for each free occurrence of $\mathrm{x}$ in
$\mathrm{A(x)}$.  The substitution is {\it free} if no free occurrence in $\mathrm{t}$
of any variable becomes bound in $\mathrm{A(t)}$; in this case we say $\mathrm{t}$
\textit{is free for} $\mathrm{x}$ \textit{in} $\mathrm{A(x)}$.\footnote{It matters which
formula was originally designated by the notation $\mathrm{A(x)}$, since if $\mathrm{x,y}$
are distinct and both occur free in $\mathrm{A(x)}$ then the sequence of free
substitutions $\mathrm{x} \mapsto \mathrm{y} \mapsto \mathrm{x}$
results in a formula different from $\mathrm{A(x)}$.}

\vskip 0.1cm

In addition to X1 - X10, {\bf Pd} has two new {\it axiom schemas}, where $\mathrm{A(x)}$
may be any formula and $\mathrm{t}$ any term free for $\mathrm{x}$ in $\mathrm{A(x)}$:
\begin{enumerate}
\item[X11.]{$\mathrm{\forall x A(x) \rightarrow A(t)}$.}
\item[X12.]{$\mathrm{A(t) \rightarrow \exists x A(x)}$.}
\end{enumerate}
In addition to R1, {\bf Pd} has two new {\it rules of inference}:
\begin{enumerate}
\item[R2.]{From $\mathrm{C \rightarrow A(x)} \;$ where $\mathrm{x}$ does not occur free
in $\mathrm{C}$, conclude $\mathrm{C \rightarrow \forall x A(x)}$.}
\item[R3.]{From $\mathrm{A(x) \rightarrow C} \;$ where $\mathrm{x}$ does not occur free
in $\mathrm{C}$, conclude $\mathrm{\exists x A(x) \rightarrow C}$.}
\end{enumerate}

\vskip 0.1cm

A {\it deduction} (or {\it derivation}) in {\bf Pd} {\it of} a formula
$\mathrm{E}$ {\it from} a collection $\Gamma$ of assumption formulas is a finite
sequence of formulas, each of which is an axiom by X1 - X12, or a member
of $\Gamma$, or follows immediately by R1, R2 or R3 from one or two
formulas occurring earlier in the sequence.  A {\em proof} is a deduction
from no assumptions.

\vskip 0.1cm

If $\Gamma$ is a collection of sentences and $\mathrm{E}$ a formula, the notation
$\Gamma \vdash_{\bf Pd} \mathrm{E}$ means that a deduction of $\mathrm{E}$ from $\Gamma$
exists.  If $\Gamma$ is a collection of formulas, we write
$\Gamma \vdash_{\bf Pd} \mathrm{E}$ only if there is a deduction of $\mathrm{E}$ from
$\Gamma$ in which neither R2 nor R3 is used with respect to any variable free in $\Gamma$.
With this restriction, the deduction theorem extends to {\bf Pd}:
If $\Gamma \cup \{\mathrm{A}\} \vdash_{\bf Pd} \mathrm{B}$ then
$\Gamma \vdash_{\bf Pd} \mathrm{(A \rightarrow B)}$.

\vskip 0.1cm

{\em Example}.  Here is a deduction in {\bf Pd} of $\mathrm{\exists x A(x)}$
from $\mathrm{\forall x A(x)}$ without using R2 or R3:
\begin{enumerate}
\item{$\mathrm{\forall x A(x) \rightarrow A(x)}$.  [axiom by X11, with $\mathrm{x}$ free
for $\mathrm{x}$ in $\mathrm{A(x)}$]}
\item{$\mathrm{\forall x A(x)}$.  [hypothesis]}
\item{$\mathrm{A(x)}$.  [by R1 from (1) and (2)]}
\item{$\mathrm{A(x) \rightarrow \exists x A(x)}$.  [axiom by X12, with $\mathrm{x}$ free
for $\mathrm{x}$ in $\mathrm{A(x)}$]}
\item{$\mathrm{\exists x A(x)}$.  [by R1 from (3) and (4)]}
\end{enumerate}
Then $\vdash_{\bf Pd} \mathrm{\forall x A(x) \rightarrow \exists x A(x)}$ follows
by the deduction theorem.

\vskip 0.07cm

Classical predicate logic {\bf Pd$^c$} comes from {\bf Pd} by strengthening
X10 to X10$^c$.  The negative interpretation extends to predicate logic using
the classical definition of $\exists$ in terms of $\forall$ and $\neg$.
The difference between constructive and classical proofs of existence is
starkly illustrated by the strong {\em existence} and {\em disjunction}
properties of {\bf Pd}:
\begin{itemize}
\item{If $\vdash_{\bf Pd} \mathrm{\exists x A(x)}$ where no variable other than
$\mathrm{x}$ is free in $\mathrm{A(x)}$, then $\vdash_{\bf Pd} \mathrm{A(x)}$ and hence
$\vdash_{\bf Pd} \mathrm{\forall x A(x)}$.}
\item{If $\vdash_{\bf Pd} \mathrm{\forall x [A(x) \vee B(x)]}$ where no variable
other than $\mathrm{x}$ is free in $\mathrm{A(x)}$ or $\mathrm{B(x)}$, then
$\vdash_{\bf Pd} \mathrm{\forall x A(x)}$ or $\vdash_{\bf Pd} \mathrm{\forall x B(x)}$.}
\end{itemize}

\subsection{Intuitionistic Predicate Logic with Equality Pd[=]}

To be useful for mathematics, the formal language must contain a binary
predicate constant $\cdot = \cdot$ denoting equality.  If $\mathrm{s, t}$ are any
terms then $\mathrm{s = t}$ is a {\em prime formula} in which all the variables
free in $\mathrm{s}$ or $\mathrm{t}$ are free.  Every prime formula of the language
of {\bf Pd} is also prime in {\bf Pd[=]}, and the {\em formulas} are built up from the
prime formulas using $\&, \vee, \rightarrow, \neg, \forall$ and $\exists$
as before.

\vskip 0.1cm

The axioms of {\bf Pd[=]} are all formulas of the extended language of the forms X1 - X12,
together with the following equality axioms, where $\mathrm{x, y}$ and $\mathrm{z}$
are distinct variables and $\mathrm{A(x)}$ may be any prime formula of the language
of {\bf Pd} in which $\mathrm{y}$ is free for $\mathrm{x}$.
\begin{enumerate}
\item[XE1.]{ $\mathrm{x = y \rightarrow (x = z \rightarrow y = z)}$.}
\item[XE2.]{ $\mathrm{x = x}$.}
\item[XE3.]{ $\mathrm{x = y \rightarrow (A(x) \rightarrow A(y))}$.}
\end{enumerate}
The rules of inference are R1 - R3 extended to the new language.
The {\em replacement property of equality} holds:  If $\mathrm{A(x)}$ is any
formula, and $\mathrm{x}$ and $\mathrm{y}$ are variables such that $\mathrm{y}$ is
free for $\mathrm{x}$ in $\mathrm{A(x)}$ and does
not occur free in $\mathrm{A(x)}$ (unless $\mathrm{y}$ is $\mathrm{x}$), then
$~\vdash_{\bf Pd[=]} ~\mathrm{x = y \rightarrow (A(x) \leftrightarrow A(y))}$.

\vskip 0.1cm

The axioms guarantee that $=$ is an equivalence relation but do
not prove that $=$ is decidable or even {\em stable} under double negation, since
$\not \vdash_{\bf Pd[=]} \mathrm{~\neg \neg (x = y) \rightarrow (x = y)}$.\footnote{Heyting
\cite{He1930B} introduced three distinct symbols to express three kinds of
equality relations:  [intensional] identity, mathematical identity, and defined
equality.  For arithmetic all three notions coincide, and number-theoretic
equality is decidable.  Equality of choice sequences, to be discussed in
the section on Brouwer's continuum, is defined extensionally and is stable
under double negation but undecidable.}

\vskip 0.1cm

A typical mathematical application will have individual constants and function
symbols, with an appropriate definition of {\em term}, and all prime formulas
will be equations.  Then the equality axioms for the function constants will
take the place of XE3, and XE2 may follow from the mathematical axioms.

\section{Intuitionistic Arithmetic HA}

Heyting \cite{He1930B} first axiomatized intuitionistic arithmetic,
which is called ``Heyting arithmetic'' in his honor.\footnote{J. van
Oosten \cite{VO2000} points out that the formalizations of {\bf HA}
as a subtheory of Peano arithmetic owe as much to G\"odel \cite{Go1933}
and to Kleene \cite{IM} as to Heyting.  Hesseling \cite{Hs2003} observes
that G\"odel used Herbrand's axioms.}  Kleene's version \cite{IM} (p. 82)
of {\bf HA} has constants and axioms for zero, successor, addition and
multiplication, and the unrestricted axiom schema of mathematical induction.

\subsection{Heyting arithmetic and Peano arithmetic}
The language of {\bf HA} is an applied version of the language of {\bf Pd[=]},
with no predicate letters, but with an individual constant $\mathrm{\bf 0}$ and function
constants $'$, $+$ and $\cdot$ .  {\em Terms} are defined inductively:
\begin{enumerate}
\item[(i)]{$\mathrm{\bf 0}$ is a {\em term}.}
\item[(ii)]{Each individual variable is a {\em term}.}
\item[(iii)]{If $\mathrm{s}$ is a {\em term}, so is $\mathrm{(s ')}$.}
\item[(iv)]{If $\mathrm{s}$ and $\mathrm{t}$ are {\em terms}, so are $\mathrm{(s + t)}$
and $\mathrm{(s \cdot t)}$.}
\item[(v)]{Nothing is a {\em term} except as required by (i) - (iv).}
\end{enumerate}
The {\em prime formulas} are the expressions of the form $~\mathrm{(s = t)}~$ where
$\mathrm{s}$ and $\mathrm{t}$ are terms.  {\em Formulas} are built up from these as usual, omitting
parentheses on the convention that $'$ binds closer than $+,\cdot$.  Every
occurrence of an individual variable in a term $\mathrm{t}$ is free in $\mathrm{t}$
(and in every prime formula containing $\mathrm{t}$).  A term or formula without
free variables is {\em closed}.

\vskip 0.1cm

The {\em axioms} of {\bf HA} are the schemas X1 - X12 for the language of
arithmetic, the equality axiom XE1, the axiom schema of {\em mathematical
induction} for arbitrary formulas $\mathrm{A(x)}$:
$$\mbox{XInd.} \hskip 0.3cm \mathrm{A({\bf 0}) ~\&~ \forall x (A(x)
\rightarrow A(x ')) \rightarrow \forall x A(x)} ,
\hskip 7.5cm$$
and the additional axioms for the primitive recursive function constants:
\begin{enumerate}
\item[XN1.]{ $\mathrm{x = y \rightarrow x ' = y '}$ .}
\item[XN2.]{ $\mathrm{x ' = y ' \rightarrow x = y}$ .}
\item[XN3.]{ $\mathrm{\neg (x ' = {\bf 0})}$ .}
\item[XN4.]{ $\mathrm{x + {\bf 0} = x}$ .}
\item[XN5.]{ $\mathrm{x + (y ') = (x + y) '}$ .}
\item[XN6.]{ $\mathrm{x \cdot {\bf 0} = {\bf 0}}$ .}
\item[XN7.]{ $\mathrm{x \cdot (y ') = (x \cdot y) + x}$ .}
\end{enumerate}

The {\em rules of inference} of {\bf HA} are R1 - R3 for the language of
arithmetic.  {\em Derivations} and {\em proofs} are defined inductively
as usual, and $\Gamma \vdash_{\bf HA} \mathrm{E}$ means that a derivation in {\bf HA}
of $\mathrm{E}$ from $\Gamma$ exists in which neither R2 nor R3 is used with respect
to a variable free in $\Gamma$.

\vskip 0.1cm

The construction of the (standard) natural numbers generates a name or
{\em numeral} for each, thus ${\bf 0}''$ is the numeral for $2$.  Each
closed term $t$ of the language expresses (under the intended interpretation)
a particular natural number; if {\bf t} is the corresponding numeral then
$\vdash_{\bf HA} \mathrm{t = {\bf t}}$.  The natural numbers are discrete:
\begin{itemize}
\item{ $~\vdash_{\bf HA} ~\mathrm{(x = y) \, \vee \, \neg (x = y)}~$.}
\item{ $~\vdash_{\bf HA} ~\mathrm{\neg \neg (x = y) \rightarrow (x = y)}~$.}
\item{ $~\vdash_{\bf HA} ~\mathrm{(x = 0) \, \vee \, \exists y (x = y')}$.}
\end{itemize}

Classical Peano arithmetic {\bf PA} comes from {\bf HA} by strengthening
X10 to X10$^c$.  G\"odel \cite{Go1933} extended the negative translation
to {\bf PA}, reducing the consistency of {\bf PA} to that of {\bf HA} and
showing that {\bf HA} cannot prove its own consistency.

\vskip 0.05cm

Kleene's proofs in \cite{IM} of G\"odel's first and second incompleteness
theorems apply to {\bf HA} as well as to {\bf PA}.  The arithmetization of
metamathematics was carried out in the intuitionistic subsystem, so
every primitive recursive predicate can be numeralwise expressed in {\bf HA}
by a {\em decidable} formula.  If $\mathrm{T(e,x,w)}$ is such a formula numeralwise
expressing ``$w$ is the g\"odel number of a computation of $\{e\}(x)$,''
and $\mathrm{A(x)}$ is $\mathrm{\exists z T(x,x,z)}$, then
$~\not \vdash_{\bf HA} \mathrm{\neg \neg \forall x (A(x) \vee \neg A(x))} $ and so
{\bf HA} +  $\mathrm{\neg \forall x (A(x) \vee \neg A(x))}$ is consistent.\footnote{To
show $\mathrm{\neg E}$ is consistent with an intuitionistic system one must show that
$\mathrm{\neg \neg E}$ (not just $\mathrm{E}$) is unprovable.}  This is a special case
of the remarkable fact that intuitionistic arithmetic is consistent with a classically
false form of Church's Thesis, as we are about to see.

\vskip 0.05cm

\subsection{Kleene's recursive realizability for intuitionistic arithmetic}

The origin and development of recursive realizability are delightfully
recounted by Jaap van Oosten in \cite{VO2000}.  Stephen Kleene, a student
of Alonzo Church, took seriously Hilbert and Bernays' assertion in \cite{HB1934}
that ``a statement of the form `there exists a number $n$ with property $A(n)$'
is \ldots an incomplete rendering of a more precisely determined proposition,
which consists either in directly giving a number $n$ with the property $A(n)$,
or providing a procedure by which such a number can be found \ldots.''  Kleene
generalized this idea to interpret every compound sentence of intuitionistic
arithmetic as an incomplete communication of an effective procedure by which
its correctness might be established, and then applied Church's Thesis to
identify ``effective'' with ``recursive.'' The result was 1945-{\em realizability}
or {number-realizability}.\footnote{Kolmogorov \cite{Ko1932} had earlier proposed
a ``problem interpretation'' of the intuitionistic connectives, but had not
connected it with recursive functions.}

\vskip 0.1cm

For the inductive definition, $n$ and $m$ range over natural numbers, and
$(n)_i$ is the exponent of the $i^{\rm th}$ prime in the complete prime
factorization of $n$ (counting $2$ as the $0^{\rm th}$ prime, and setting
$(0)_i = 0$ by convention.  Ordered pairs are coded by $(n,m) = 2^n \cdot 3^m$,
and $\{n\}(m)$ denotes the result of applying the $n^{\rm th}$ recursive partial
function to the argument $m$.

\vskip 0.1cm

{\em Definition}. (Kleene 1945)  A number $n$ {\em realizes} a sentence $\mathrm{E}$
only as follows:
\begin{enumerate}
\item{$n$ {\em realizes} a closed prime formula $\mathrm{r = t}$, if $\mathrm{r = t}$
is true under the intended interpretation.}
\item{$n$ {\em realizes} $\mathrm{A ~\&~ B}$, if $(n)_0$ {\em realizes} $\mathrm{A}$
and $(n)_1$ realizes $\mathrm{B}$.}
\item{$n$ {\em realizes} $\mathrm{A \vee B}$, if either $(n)_0 = 0$ and $(n)_1$
{\em realizes} $\mathrm{A}$, or $(n)_0 \neq 0$ and $(n)_1$ {\em realizes} $\mathrm{B}$.}
\item{$n$ {\em realizes} $\mathrm{A \rightarrow B}$, if, for every $m$:  if $m$
{\em realizes} $\mathrm{A}$ then $\{n\}(m)$ is defined and {\em realizes} $\mathrm{B}$.}
\item{$n$ {\em realizes} $\mathrm{\neg A}$, if no $m$ {\em realizes} $\mathrm{A}$.}
\item{$n$ {\em realizes} $\mathrm{\forall x A(x)}$, if, for every $m$: $\{n\}(m)$ is
defined and {\em realizes} $\mathrm{A({\bf m})}$.}
\item{$n$ {\em realizes} $\mathrm{\exists x A(x)}$, if $(n)_1$ {\em realizes}
$\mathrm{A({\bf m})}$ where $m = (n)_0$.}
\end{enumerate}

A formula $\mathrm{E}$ is {\em realizable} if some $n$ realizes the universal closure
$\mathrm{\forall E}$ of $\mathrm{E}$.

\vskip 0.1cm

Kleene conjectured that every closed theorem of {\bf HA} was realizable.
His student David Nelson verified the conjecture, then formalized the proof
in an extension of {\bf HA}.\footnote{The technical work of formalization,
in \cite{Ne1947}, was nontrivial.  Kleene \cite{Kl1945} announced and
interpreted Nelson's results.  For a comprehensive and comprehensible
modern treatment, see Troelstra \cite{Tr1998}.}  The key lemma states that
for each closed term $\mathrm{t}$ expressing the number corresponding to the
numeral $\mathrm{\bf t}$:
$$\mbox{$n$ realizes $\mathrm{A(t)}$  } \Leftrightarrow
\mbox{  $n$ realizes $\mathrm{A({\bf t})}$.}$$

\vskip 0.1cm

{\em Nelson's Theorem}.   If
$\mathrm{C}_1, \dots, \mathrm{C}_k \vdash_{\bf HA} \mathrm{A}$ and
$\mathrm{C}_1, \ldots, \mathrm{C}_k$ are realizable, so is $\mathrm{A}$.

\vskip 0.1cm

{\em Corollary 1}.  If $~\vdash_{\bf HA} \mathrm{\forall x}_1 \ldots \mathrm{\forall x}_n
\mathrm{\exists y A(x}_1,\ldots,\mathrm{x}_n,\mathrm{y)}~$ where
$~ \mathrm{A(x}_1,\ldots,\mathrm{x}_n,\mathrm{y)}~$ contains free only
$\mathrm{x}_1,\ldots,\mathrm{x}_n,\mathrm{y}$ then there is a general recursive function
$\psi$ of $n$ variables such that for all values of $\mathrm{x}_1,\ldots,\mathrm{x}_n$:
If $~\psi(x_1,\ldots,x_n) = y$, then
$\mathrm{A({\bf x}}_1,\ldots,\mathrm{\bf x}_n,\mathrm{{\bf y})}$
is realizable.

\vskip 0.1cm

In \cite{Kl1945} Kleene observed that the cases of the definition for
$\vee$, $\rightarrow$ and $\exists$ could be modified to give another notion
(later called {\em realizable-\,$\vdash$}) for which the analogue of Nelson's
Theorem held, with the following result.\footnote{The existence and disjunction
properties for {\bf HA} were implicit special cases, as Kleene later noted.}

\vskip 0.1cm

{\em Corollary 2} (to the version of Nelson's Theorem for ``realizable-\,$\vdash$'').
Assuming $~\vdash_{\bf HA} \mathrm{\forall x}_1 \ldots \mathrm{\forall x}_n
\mathrm{\exists y A(x}_1,\ldots,\mathrm{x_n,y)}~$
where $~ \mathrm{A(x}_1,\ldots,\mathrm{x}_n,\mathrm{y)}~$ contains free only
$\mathrm{x}_1,\ldots,\mathrm{x}_n,\mathrm{y}$,
there is a general recursive function $\psi$ such that for all values of
$x_1,\ldots,x_n$:
$$~\vdash_{\bf HA} \mathrm{A({\bf x}}_1,\ldots,\mathrm{\bf x}_n,\mathrm{{\bf y})}~ \;
\mbox{where} \; ~\psi(x_1,\ldots,x_n) = y.$$

\subsection{Church's Thesis}

It is possible to express Church's Thesis in the language of arithmetic.
One version (which includes countable choice) is the schema CT$_0$:
$$\mathrm{\forall x \exists y A(x,y) \rightarrow
\exists e \forall x \exists w [T(e,x,w)~\&~A(x,U(w))]}$$
where $\mathrm{T(e,x,w)}$ (numeralwise) expresses ``$w$ is the g\"odel number of a
computation of $\{e\}(x)$,'' and $\mathrm{A(x,U(w))}$ abbreviates
$\mathrm{\forall z (U(w,z) \rightarrow A(x,z))}$ where $\mathrm{U(w,z)}$ (numeralwise)
expresses ``$z$ is the value, if any, computed by the computation with g\"odel number $w$;
otherwise $z = 0$''.  The g\"odel numbering is primitive recursive, and
$\mathrm{T(e,x,w)}$ and $\mathrm{U(w,z)}$ are quantifier-free.\footnote{For
details including the definition of {\em numeralwise expressibility} see \cite{IM}.}
Nelson's Theorem entails the following consistency and independence results.

\vskip 0.1cm

{\em Corollary 3}.  {\bf HA} + CT$_0$ is consistent.

\vskip 0.1cm

{\em Corollary 4}.  If $\mathrm{A(x)}$ is $\mathrm{\exists y T(x,x,y)}$ and
$\mathrm{B(x)}$ is $\mathrm{A(x) \vee \neg A(x)}$, then
\begin{enumerate}
\item[(i)]{$~\not \vdash_{\bf HA} \mathrm{\forall x (A(x)~\vee~\neg A(x))}$,}
\item[(ii)]{$~\not \vdash_{\bf HA} \mathrm{\neg \neg \forall x (A(x)~\vee~\neg A(x))},\;$
and}
\item[(iii]{$~\not \vdash_{\bf HA} \mathrm{\forall x \neg \neg B(x) \rightarrow
\neg \neg \forall x B(x)}$.}
\end{enumerate}

While {\bf HA} is a subsystem of {\bf PA}, {\bf HA} + CT$_0$ is a nonclassical
arithmetic.  To see why, let $\mathrm{A(x,y)}$ be
$\mathrm{(y = 0 \rightarrow \forall z \neg T(x,x,z))~\&~(y \neq 0 \rightarrow
T(x,x,y\dot{-}1))}$.  Evidently $\vdash_{\bf PA} \mathrm{\forall x \exists y A(x,y)}$.
In {\bf HA} the hypothesis $\mathrm{\forall x \exists w [T(e,x,w)~\&~A(x,U(w))]}$
implies $\mathrm{\exists w [T(e,e,w)~\&~A(e,U(w)]}$, but
$\vdash_{\bf HA} \mathrm{\forall w [T(e,e,w) \rightarrow A(e,w+1)]}$ and the g\"odel
numbering satisfies $~\vdash_{\bf HA} \mathrm{\forall w (U(w) \leq w)}$ and
$\vdash_{\bf HA} \mathrm{\forall u \forall v [A(e,u)~\&~A(e,v) \rightarrow u = v]}$.
Hence
$~\vdash_{\bf HA} \mathrm{\neg \exists e \forall x \exists w [T(e,x,w)~\&~A(x,U(w))]}$,
so {\bf PA} proves the negation of an instance of CT$_0$.

\subsection{Axiomatization and modifications}

When formalizing the original definition, Nelson associated with each
formula $A$ of {\bf HA} another formula $\mathrm{e}$ {\sc r} $\mathrm{A}$ of (a
conservative extension of) {\bf HA}, then proved that every formula of the form
$\mathrm{A \leftrightarrow \exists e (e}$ {\sc r} $\mathrm{A)}$ was realizable and hence
consistent with {\bf HA}.  In 1971 Troelstra used an extension ECT$_0$ of
Church's Thesis CT$_0$ to determine the exact strength of number-realizability
over intuitionistic arithmetic.\footnote{ECT$_0$ is the schema
$\mathrm{\forall x [A(x) \rightarrow \exists y B(x,y)] \rightarrow
\exists e \forall x [A(x) \rightarrow \exists w (T(e,x,w) ~\&~ B(x,U(w)))]}$,
where $\mathrm{A(x)}$ is {\em almost negative} (containing no $\vee$, and no $\exists$
except immediately in front of a prime or quantifier-free formula); every formula
$\mathrm{e}$ {\sc r} $\mathrm{A}$ is of this kind.}  Briefly, Troelstra showed that
the provable sentences of {\bf HA} + ECT$_0$ are exactly those whose realizability
can be established in {\bf HA}, and that
$\vdash_{{\bf HA} + {\rm ECT}_0} \mathrm{(E \leftrightarrow \exists x (x}$ {\sc r}
$\mathrm{E))}$ for every formula $\mathrm{E}$.

\vskip 0.1cm

From the intuitionistic point of view Church's Thesis is restrictive, and
probably unacceptable as a general principle.  Markov's Principle MP, the schema
$$\mathrm{\forall x (A(x) \vee\neg A(x)) \rightarrow [\neg \forall x \neg A(x)
\rightarrow \exists x A(x)],}$$ is also problematic.  Because MP implies its own
realizability, Troelstra's theorem extends over {\bf HA} + MP to
{\bf HA} + ECT$_0$ + MP.\footnote{Troelstra has convincingly argued that the
``Russian recursive mathematics'' of the Markov school is based on this theory.}
Kreisel \cite{Kr1958} invented a typed version of realizability in order to show
the independence of Markov's Principle.

\vskip 0.1cm

Other modifications of realizability have been developed to prove a variety
of independence and consistency results, and this process continues.
The original notion gives a classically comprehensible interpretation of
Heyting arithmetic but does not claim to capture intuitionistic arithmetical truth.

\vskip 0.5cm

\section{The intuitionistic theory of the continuum}

Brouwer's main objection to classical mathematics (apart from the unrestricted
use of the principle of excluded middle) was ``its introduction and description
of the continuum.''\footnote{Changes in the Relation between Classical Logic
and Mathematics, in \cite{VS1990} p. 453.}  His entire work was motivated by
the attempt to describe a construction of the continuum in harmony with his
mathematical principles, and to develop a satisfactory mathematics on the basis
of that construction.\footnote{As Beth says (\cite{Be1959} p. 422), ``the
central theme in intutionistic mathematics is the theory of the continuum.''}
To achieve this goal he created the new notions which give intuitionistic mathematics
its unique character.

\vskip 0.1cm

Our idea of the continuum, the real plane or in the one-dimensional case the real line,
seems to originate from our perception of space, a primary notion
which is {\em a priori} according to Kant.  Thus, until non-Euclidean
geometries appeared and the uniqueness of Euclidean space was lost, the
continuum could be considered an initial concept, immediately comprehensible
by intuition, not requiring analysis in terms of other more elementary concepts.
After this loss, however, mathematicians began trying to define the continuum
from apparently more fundamental notions, by more and more abstract methods.
It was the era of Cantor's creation of the theory of sets, with the ordinal
and cardinal numbers, and of the arithmetization of analysis by Weierstrass,
Dedekind, Cantor and others.

\vskip 0.1cm

Since then, in traditional classical mathematics the continuum is considered
to be a collection of distinct mathematical objects, the real numbers, which
are defined using Dedekind cuts or Cauchy sequences of rational numbers, or nested
sequences of intervals with rational endpoints.  In each case, infinite
mathematical entities are treated as completed, actual.  Even the (semi)intuitionists
Poincar\'e, Borel and Lebesgue, in their treatment of the continuum, failed to
maintain the constructivist standpoint.  Brouwer criticized them for having
``recourse to logical axioms of existence'' such as the axiom of completeness,
or ``[contenting] themselves with an ever-denumerable and ever-unfinished''
(\cite{Br1981} p. 5) set of numbers so that, in order to retain constructivity,
they were restricted to definable and hence to only countably many
objects (for example, definable sequences of rationals).\footnote{In addition
to the problem of cardinality, the measure of a countable set of reals is zero,
and that problem also concerned Brouwer from the beginning.}

\vskip 0.1cm

Brouwer accepted as correct Cantor's argument showing that it is
impossible to enumerate all the points of the continuum.  His first answer
to the resulting foundational problem was presented in his dissertation, where
he considered the continuum as a primitive entity, ``the inseparable complement
of the discrete,'' to which it cannot be reduced:
\begin{quote}
``However, the continuum as a whole was given to us by intuition; a construction
for it, an action which would create from the mathematical intuition `all' its
points as individuals, is inconceivable and impossible'' because the
``mathematical intuition is unable to create other than denumerable sets of
individuals.'' (\cite{Br1907} p. 45)
\end{quote}
According to its intuitive conception, the continuum can be infinitely divided:
two discrete points are connected by a ``between", which is never exhausted
by the insertion of new points.  So it is possible ``after having created
a scale of ordertype ç [the ordertype of the rationals] to superimpose upon it
a continuum as a whole" (\cite{Br1907}, p. 45).  Of course, this description of
the manner in which we appear to pass from the rationals with their dense linear order
to the continuum does not give a construction of the continuum, it simply indicates
the relationship between the rationals and the continuum and perhaps differs little
from the acceptance of an axiom of completeness, since the gaps between the rational
points are covered by their (rather indefinite) ``between.''
In any case the intuitionistic construction of the continuum remained a challenge
for Brouwer; his answer to that problem was given later and stands as his
most important creation.

Among various thoughts concerning the continuum, in his
doctoral thesis Brouwer presented arguments for the continuum hypothesis, i.e. the
conjecture that Cantor's second number class and the continuum have the same number
of points:  this holds primarily, in Brouwer's early opinion, because considering the
continuum as a set of points we are obliged to define its points somehow, and therefore
a set of the same cardinality as the second number class results.\footnote{Indeed
Brouwer notes that the problem has not been well posed, because ``neither the second
number class nor the continuum as a totality of individualized points exists
mathematically'' (\cite{Br1907}, p. 83).}

\vskip 0.1cm

Thus, in Brouwer's thought as expressed in his dissertation, the continuum
had an intuitive basis very close to its geometrical representation,
like that of the semi-intuitionists.  A constructive \textit{arithmetical}
(in the sense of being based in some way on some notion of number) view of
the continuum was impossible, and so its mathematical treatment remained
problematic.

In his lectures of his mature period, Brouwer gave to his
first, essentially restrictive intervention in the problems of foundations
of mathematics the name ``First Act of Intuitionism."\footnote{\textit{The
First Act of Intuitionism.}  The complete separation of mathematics from
mathematical language, and hence from the linguistic phenomena which are
described by theoretical logic, the recognition that intuitionistic mathematics
is an essentially languageless activity of the mind which springs from the
perception of a motion of time.  This perception of time can be described as
the splitting of a life moment into two distinct things, one of which gives
place to the other, but is preserved in memory.  If the dyad thus born is
stripped of every quality, what remains is the blank form of the common
substratum of all dyads.  And it is this common substratum, this common form,
which is the basic intuition of mathematics. (\cite{Br1981}, p. 4)}
With it he completely separated language and logic from mathematics and
ruled that the intuition of many-oneness is the only basis for
mathematical activity - as he had already stated in his dissertation.
He recognized however the restricted possibilities of mathematical development
(only ``separable'' mathematics like arithmetic and
algebra survived):  ``Since the continuum appears to remain outside
its scope, one might fear at this stage that in intuitionism would be no place
for analysis.'' (\cite{Br1981}, p. 7).  And then he presented
the ``\textit{Second Act of Intuitionism},''\footnote{\textit{The Second Act of
Intuitionism.} The adoption of two ways of creating new mathematical entities:
first in the form of infinite sequences of mathematical entities which have
already been constructed, which continue more or less freely (so that, for
example, infinite decimal fractions which neither have exact values nor are
guaranteed to ever acquire exact values are allowed); second, in the form of
mathematical species, that is properties which are supposed to hold of mathematical
entities which have already been constructed, which satisfy the condition that
if they hold for a mathematical entity, they also hold for all the mathematical
entities which have been determined to be `equal' to it, where the definitions
of equality must satisfy the conditions of symmetry, reflexivity
and transitivity. (\cite{Br1981}, p. 8)} with which he introduced
his new notions of \textit{free choice sequences} and \textit{species} and
began the intuitionistic reconstruction of analysis.

\subsection{Brouwerian set theory}

\subsubsection{Choice sequences}
 As we have seen, the role of infinite sequences is crucial in an arithmetical
description of the continuum.  Borel had remarked that the only way to obtain
the continuum using only sequences of rationals without imposing the existence
of real numbers by means of axioms, was to adopt sequences of arbitrary choices
of objects.  As a constructivist he hesitated but did not completely reject the
idea:
\begin{quote}
``...one knows that the completely arithmetical concept of the continuum
requires that one admits the legitimacy of a countable infinity of successive
choices.  This legitimacy seems to me highly debatable, but
nevertheless one should distinguish between this legitimacy and the legitimacy
of an uncountable infinity....The latter concept seems to me ... entirely
meaningless..." (\cite{Bo1909} and \cite{TvD1988} volume 2, p. 641).
\end{quote}
Brouwer, with his Second Act of Intuitionism, accepted free choice sequences
as legitimate mathematical objects.\footnote{This acceptance reflected his
opinion, inspired by Kant, that mathematical constructions are internal,
free creations of the mind of the \textit{creating subject}, a view which
influenced his methodology.}  As we shall see, he found a way to use these
infinite, undetermined objects constructively by viewing them as having two parts:
a finite, already constructed part which permits genuine constructive use of the
infinite sequence under some circumstances, and an infinite, undetermined part which
makes it possible to obtain the whole continuum, escaping the restrictions
imposed by any sort of definability.

\vskip 0.1cm

\subsubsection{The new notions of set.}  Brouwer's new perspective on the
foundations of mathematics was presented in his 1918 article ``\textit{Foundation
of Set Theory, independent of the logical law of excluded middle.  First part:
General Set Theory}.'' \cite{Br1918}.  There, instead of Cantor's sets,
two alternative notions are proposed as the basis for analysis, the \textit{spread}
(Menge)\footnote{Brouwer initially used the term Menge, which had been preempted
by Cantor, but later preferred the term spread.  In the Greek version of this
article we translated ``spread'' using an ancient Greek word for loom, following
the example of Wim Veldman (possibly the unique modern mathematician who works
precisely in the spirit of Brouwer).} and \textit{species}.

\vskip 0.1cm

\textit{Species:  the set as property.}  Species are in fact closely related to
sets as understood in classical mathematics.  A species is a property, but a
property ``supposable for mathematical entities previously acquired'', and is
extensional with respect to the notion of equality (in general an
equivalence relation) between the mathematical objects in question.  The verification
that a mathematical object has a property (hence ``belongs to the species'')
requires a construction; the species therefore may be thought of as a construction
\textit{within} a construction.  Two species are the same ``if to every element
of each an equal element of the other can be correlated.''  With this definition
of species, the problems of impredicativity and self reference are avoided:
``a species can be an element of another species, but never an element of itself!''
as Brouwer himself asserts (\cite{Br1981}, p. 8).

\vskip 0.1cm

\textit{Spread:  the set as law.}  But apart from the natural numbers, which
mathematical entities can be elements of species?  The generation of mathematical
objects results from mathematical activity typically involving
the radically new concept of spread, where free choice sequences are allowed.
Each element of a spread is constructed in stages: at each stage
we choose, in accordance with certain rules, one more piece of the element, thus
obtaining a better approximation of it.  Specifically:

\vskip 0.1cm

A spread is created i) by successive choices of natural numbers;
each choice depends on the preceding ones, and is made either freely
or in accordance with certain restrictions, and ii) by assigning after each
choice an object from a particular preexisting denumerable set.  The
restrictions in i) and the correlations in ii) are given by the
\textit{spreadlaw}, in reality a pair of laws, the \textit{choice law}
(which also determines whether the process will terminate or be continued, and
in the second case there must be at least one permitted choice for the
next step)\footnote{This means that we can assume that only infinite
sequences are produced by the spread law.} and the \textit{correlation law}.
Attempting to describe the nature of the `laws', Brouwer notes:
``...we can say that a spread law yields an \textit{instruction}\footnote{Italics
ours.} according to which ..."
(\cite{Br1981}, p. 15).

\vskip 0.1cm

The \textit{elements} of the spread are the sequences (infinite or
finite) of the objects correlated to each sequence of choices.
So an element may remain forever incomplete, always becoming.
The spread thus created is not considered to be the
totality of its elements, but is identified with the law of its creation.
The description of the notion of spread is complicated and long, but Brouwer
said he could not avoid this complexity.  Over time he made various modifications
to the details of the definition, concerning for example the nature of
the restrictions imposed.

\vskip 0.1cm

Despite the fact that every choice and subsequent assignment is made
effectively and so it is decidable whether a given finite sequences of choices
or of objects to be assigned is permitted, the same is not true for the
infinite sequences: a free choice sequence is constructed in steps; if at some
step the choice law determines that the finite part which has already been
created no longer satisfies the restrictions, then indeed we know that the
sequence does not `belong' to the spread.  But as long as it satisfies the
restrictions, we can never be sure that it will not be rejected at some later
step, and that is the price of freedom in the creation of the sequence.

\vskip 0.1cm

Spreads are more basic objects than species: the spread generates its elements,
while the species has elements already constructed, and what is needed for a
mathematical entity to belong to a species is the intuitionistic justification
that it has the property which defines the species.

\subsubsection{Spreads as trees with a topology.}  We can think of spreads as
trees, whose nodes are the finite sequences of correlated objects corresponding
to the successive choices, and all of whose branches are infinite (a branch is a
sequence each of whose finite initial segments belongs to the tree).  Brouwer
himself used this image in describing spreads.

\vskip 0.1cm

In a very natural way a topology can be defined on a spread:  considering the
infinite branches (using the tree image) as points and taking
as basis for the topology the set $\mathcal{N}_u$, where $\mathcal{N}_u$ is
the set of the infinite branches with initial segment $u$.

\subsubsection{Examples.}  Now we give two most important examples of spreads.

\vskip 0.1cm

 1.  The \textit{universal spread}, whose choice law
permits the choice of any natural number whatsoever at each step, and whose
correlation law assigns trivially after each choice the number just
chosen.  This spread gives us all the sequences of natural numbers, so the
species of its elements is uncountable.  Equipped with the initial segment
topology, it becomes the familiar Baire space, and so it is Borel-isomorphic
(isomorphic as a measure space) to the set $\mathcal{R}$ of the traditional
real numbers; it is (topologically) homeomorphic with the irrational numbers.
This most general and free process associated with the universal spread captures
the essence of the creation of the intuitionistic continuum, and on this basis,
by placing a few minor restrictions, we can define the real numbers in the
context of intuitionism.

\vskip 0.1cm

2. A finitely branching spread is called a \textit{fan}.  Fans, which are
compact in the above topology, play a special role in the development of a theory
of analysis because of their more determinate character.  One fan of particular
significance is the binary fan, consisting of all infinite sequences of 0s
and 1s, which with the initial segment topology is the well-known Cantor space.

\vskip 0.1cm

The intuitionistic continum and the real numbers are the next examples, to which
we now turn.

\subsection{The intuitionistic continuum.}
\subsubsection{The spread of the points of the continuum and the real numbers.}
Let us consider\footnote{From the footnote of Parsons in \cite{Br1927}.} the
species of the binary fractions $a/2^n$, where $a$ is an integer and $n$ a
natural number, with the usual ordering.  In this species, we consider the closed
intervals of the form $I_{m,n}=[m/2^n,(m+2)/2^n]$ and an enumeration of the pairs of
natural numbers $p_1, p_2,...$ with $p_i=<r_i,s_i>$, and the following spread: at
the first step the choice of any natural number $n$ is permitted and the interval
$I_{r_n,s_n}$ is assigned; if $a_1,...,a_n$ have been chosen and the interval
$I_{r,s}$ has been assigned, then it is permissible to choose the number $k$ only
if $I_{r_k,s_k}$ is properly contained in $I_{r,s}$, and if $k$ is chosen then
the interval $I_{r_k,s_k}$ is assigned.  The elements of this spread, these
infinitely proceeding sequences of closed intervals whose lengths converge to 0,
are \textit{the points of the continuum}.\footnote{We can compare this process
with the \textit{nested interval principle}, a completeness axiom which assures
that the intersection of every sequence of nonempty closed nested intervals whose
lengths converge to 0 is nonempty.}

\vskip 0.1cm

An equivalence relation is now defined for these points:  two points $p$
and $q$ \textit{coincide} if every interval of $p$ intersects
every interval of $q$.  Each equivalence class is a species called a
\textit{point core} and this is what a real number is defined to be.  The
species of all these point cores is \textit{the intuitionistic continuum}.
The equality between reals is obtained from the coincidence relation:  if $r$ and $s$
are reals and $p, q$ representatives of $r, s$ correspondingly, then $r = s$ if
and only if $p$ and $q$ coincide.

Brouwer gave many analogous constructions of the continuum, using for example sequences
of rationals or binary fractions with various rates of convergence.  Heyting
also gave some constructions of this kind, using the characteristic term
\textit{real number generators} for the corresponding point cores.  This
description of the real numbers solves the problems of cardinality and measure
from an intuitionistic standpoint.  For example, the rationals are naturally embedded
into the reals defined as above, since every rational belongs to a
sequence of intervals $I_{m,n}$ with lengths tending to $0$. But the differences from
the real numbers of classical mathematics remain large, as will appear in what follows.

\subsubsection{Undecidability of equality.}  In his 1930 paper titled
``\textit{The structure of the continuum}'' \cite{Br1930} Brouwer examined
basic properties of the continuum, comparing the classical continuum with
the one resulting from his own views.  His first conclusion was
that equality between two real numbers \textit{is undecidable} in the case
of the intuitionistic continuum.  We can see this as follows:  Let $A(n)$
be a property which is decidable for each $n$, for example the Goldbach
conjecture (every even number greater than 3 is the sum of two primes),
so that $A(n)$ is the property that $2n+4$ is the sum of two primes. Let
$$\alpha (n)=\left\{
\begin{array}{ll}
1/2^n, \;\;&\mathrm{if} \; \;\;\forall k\leq n \;A(k), \\
1/2^k, &\mathrm{if} \; \;\neg A(k)~\&~ k \leq n ~\&~ \forall m<k
\;A(m).
\end{array}\right.$$
The sequence of rationals so defined is evidently convergent, hence determines
a real number  $r$.
We see that $r = 0$ if and only if $A(n)$ holds for every $n$
(otherwise $r = 1/2^k$ for the least $k$ for which $A(k)$ fails).  So we
cannot decide if $r$ coincides with $0$, as long as we do not know the answer to
Goldbach's conjecture.

\subsubsection{Weak counterexamples.}  Counterexamples of the type just
considered are characteristic of Brouwer's argumentation, and are known as
\textit{weak counterexamples}.  Each one uses a problem which has not yet been
solved in order to deny a classically valid statement.  It is weak because it
depends crucially on a \textit{specific} unsolved problem, to which a solution
may someday be found; for instance, many of Brouwer's weak counterexamples
involved Fermat's Last Theorem, which was then undecided but has now been proved.

But Brouwer believed that there would always be unsolved mathematical problems.
A weak counterexample can be considered as a construction which will convert
any unsolved problem into a refutation of a classical theorem.  Implicit in
the use of weak counterexamples is Brouwer's continuity principle, to be
explained below.

\subsubsection{The continuum cannot be ordered.} A second result was the
impossibility of finding a linear order of the intuitionistic continuum.
The rationals are obviously totally ordered, because the comparison of
rationals reduces to the comparison of natural numbers.  For the reals,
some can certainly be compared and ordered among themselves;
in the characterization by nested intervals for example, if at some
step of the construction of two points of the continuum $p$ and $q$
an interval of $p$ lies entirely to the left of an interval of $q$,
then, for the corresponding reals $r$ and $s$ we have $r<s$. Thus
the \textit{natural order} of the reals, as Brouwer called it, is
derived. But again by appealing to weak counterexamples, he showed
that this order cannot be total. Then he defined a refinement of the
natural order, the \textit{virtual order}
     $$r\prec s\equiv \neg r > s ~\&~ r\neq s, $$
which also fails to be total.  It is worth noting that,
although from the intuitionistic standpoint not all real numbers are comparable,
the natural ordering can be proved to satisfy certain properties;
the most useful substitute for comparability is the property
$r<s\rightarrow r<t \vee t<s$.  The lack of order is not a ``technical''
problem, but rather a consequence of the fact that most real numbers are
incompletely determined objects.

\subsubsection{The unit continuum.}  Because two real numbers $r$ and $s$
are not always comparable with respect to the natural ordering, the closed
interval $[r,s]$ is defined to be the species of all reals $x$ for which
it is impossible that $x > r$ and $x > s$, and also impossible that $x < r$
and $x < s$.  If $r < s$ then $[r,s]$ is just the species of $x$ such that
neither $x < r$ nor $x > s$.
The proof of the uniform continuity theorem, discussed below, depends on the fact
that any closed interval $[r,s]$ is the continuous image of a fan F of real number
generators, each of which coincides with a point of $[r,s]$.  By the {\em unit continuum}
Brouwer meant a fan coordinated in this way with $[0,1]$ (as in \cite{Br1981} p. 35;
cf. \cite{TvD1988} and \cite{Lo2005}).


\subsection{The basic theorems and implicit principles}

\subsubsection{Discontinuous functions do not exist.}  The undecidability
of equality has consequences which a traditional mathematician would find
hard to accept.  Let\footnote{The example is from \cite{TvD1988}, volume 1,
p. 14.} $f$ be the function from the reals to the reals with
$$ f(x)=\left\{\begin{array}{ll}
1\;\;&\mathrm{if}\;\;x\neq 0,\\
0\;\;&\tx{otherwise,}
\end{array}\right.$$
which is classically defined everywhere and discontinuous at 0.  Intuitionistically
however, $f$ cannot be a total function: for the real number $r$ of the weak
counterexample we gave, as for any other real, it follows by continuity
(in the usual sense) that we can calculate any approximation to $f(r)$, and from
the presumed totality of the function we can decide whether $f(r) < 1$ or
$f(r) > 0$, and so decide whether   $\neg r \neq 0$ or $r \neq 0$, or
equivalently\footnote{The equivalence of   $\neg r \neq 0$ with $r = 0$ needs a little
argument.  The equality $r = 0$ is of the form $\forall n \;A(n)$, and the $A(n)$
is decidable hence stable under double negation, while the implication
$\neg \neg \forall n A(n) \rightarrow \forall n \neg \neg A(n)$ holds
intuitionistically.}
whether $r = 0$ or $r \neq 0$, which is impossible; thus we arrive at the conclusion
that $f$ cannot be a total function.

\vskip 0.1cm

Another impressive result obtained by similar arguments is that
the Intermediate Value Theorem fails intuitionistically.

\subsubsection{The uniform continuity theorem.} In contrast with these negative
results, Brouwer reached a very strong conclusion, similar to the classical
theorem that functions continuous on a compact space are uniformly continuous,
except that in the intuitionistic case only continuous functions exist:

\vskip 0.1cm

\textit{Uniform continuity theorem}.  Every total function defined on the
unit continuum is uniformly continuous.

\vskip 0.1cm

In fact Brouwer probably arrived at this result guided by his intuition
rather than by clear mathematical arguments.  Hermann Weyl said the following
in one of his lectures, defending Brouwer's continuum:  ``It is clear that
no one can explain the meaning of `continuous function on a bounded interval'
without including `uniform continuity' and boundedness in
the definition.  Above all, there can be no function on a continuum other
than continuous functions.'' (cite{VS1990}, p. 379).  Brouwer read and
emphatically agreed with these remarks.  For many years he tried to find
a satisfactory proof for this theorem.  In his 1927 article ``\textit{On the
domains of definition of functions}'' \cite{Br1927} the most complete
presentation of his arguments is given.  The value of this attempt
consists mainly in the fact that in his proof, especially in its final form,
the two characteristic principles of intuitionistic analysis (which we will
discuss below) are used clearly and thus brought to light, though not explicitly
as principles but only as natural manifestations of the intitionistic viewpoint.

\subsubsection{The Fan Theorem.}  The uniform continuity theorem was obtained
as a consequence of the fan theorem, as Brouwer called it.

\vskip 0.1cm

\textit{The Fan Theorem}.  If to each element $\alpha$ of a fan a number
$b_\alpha$ is associated, then a natural number $z$ can be found so that
$b_\alpha$ is determined completely by the first $z$ values of $\alpha$.

\vskip 0.1cm

This version of the theorem fails classically; its proof uses the
(classically correct) method of backwards induction but also the (classically
false) intuitionistic continuity principle\footnote{The principles mentioned here
will be defined below.}.  It is the form which Brouwer preferred.

\vskip 0.1cm

A different version of the theorem is provable just using the principle of
backwards induction and is classically correct.  In this version, the theorem
says that if, for the given fan there is a decidable set $B$ of nodes such that each
branch $\alpha$ meets $B$ (such a $B$ is called a \textit{bar}), then there is
a number $z$ such that each branch meets the bar at a node of length at
most $z$.  We mention this alternative because it is classically equivalent to
K\"onig's Lemma (every infinite, finitely branching tree has an infinite branch).

\vskip 0.1cm

As van Dalen notes in \cite{Dal1999}, ``it is an interesting historical curiosity
that the fan theorem preceded its much better known contrapositive:
K\"onig's infinitary lemma [K\"onig 1926] ... The infinitary lemma is not
constructively valid.  An interesting observation is that this contradicts the
popular impression that constructive proofs and theorems are always later
improvements of classical theorems.''

\subsubsection{The Bar Theorem.}  In order to prove the fan theorem Brouwer
first gave a complicated proof, more or less metamathematical in character, of
what he later called the \textit{bar theorem}, in which he examined the structure
of possible proofs of sentences of certain forms.  Later however, as appears in a
footnote to his 1927 article \cite{Br1927}, he realized that what he uses
is in fact a kind of induction principle.  A statement of this principle
follows:

\vskip 0.1cm

\textit{Bar Theorem}.  Let the universal spread contain a decidable bar $A$
and $X$ be a set of nodes satisfying
(i) $A \subset X$ and (ii) if all the immediate successors of a node $n$
belong to $X$, then $n$ belongs to $X$; then it follows that the empty
sequence (the root of the universal spread) belongs to $X$.

\vskip 0.1cm

We remark that this principle makes it possible to exploit properties possessed
by all sequences of natural numbers, despite their incomplete character,
in the case where these properties can be verified at some finite stage
of the generating process of the sequence.  But in any case the assumption of
the existence of a bar on the universal spread is not at all trivial; to see
this, we need only observe that the set of nodes where each branch meets
the bar for the first time may be as ``long'' as any infinite countable
ordinal number.  While this principle seems to be closely connected with
intuitionistic mathematics, it also happens to be classically correct.

\subsubsection{The continuity principle.}  When Brouwer proved that the
universal spread generates uncountably many elements, he replaced Cantor's
diagonal argument by a new argument of exceptional simplicity, as follows:
suppose we have a function from the universal spread to the natural numbers,
and let $b$ be the value this function takes on the sequence $\alpha$.
Constructively, this $b$ must be determined by the first $y$ values of $\alpha$,
for some $y$.  But then the function must also take the value $b$ on a sequence
$\beta$ having the same first $y$ terms as $\alpha$, although $\beta(y+1) \neq
\alpha(y+1)$.  Hence it is impossible to map the universal spread to the
natural numbers in a one-to-one way, and so the elements of the universal
spread cannot be enumerated.

\vskip 0.1cm

The substance of this argument is \textit{the continuity principle},
which played a crucial role in the proofs of several theorems considered above.
According to this principle, in order for a function to be defined on a spread,
the value it takes on each sequence-argument must be determined entirely
by an initial segment of the sequence.  This is how the conflict
between the incompletely determined nature of choice sequences, and the
constructive character of a function defined over them, was reconciled by
Brouwer.  So the principle postulates that every total function is continuous.
It can also be generalized to functions with sequence values.  It is the
only genuinely intuitionistic principle of Brouwer's mathematics, it
contradicts classical analysis, and thus the mathematics based on it
are distinct from, in fact inconsistent with, classical mathematics.

\subsubsection{Other mathematical principles.}  It is certain that Brouwer
saw the general axiom of choice as eminently nonconstructive.  However various
arithmetical forms of it result by correct reasoning from the manner in which
existential and universal statements are understood intuitionistically.
We will discuss these in the context of a formal system which was proposed
for Brouwer's analysis.

\subsection{$\mathcal{FIM}$:  a formal system for intuitionistic analysis}
The axiomatization of intuitionistic mathematics by Heyting in his classic
work \cite{He1930A, He1930B, He1930C} included the theory of sets, but in a
form which, as Kleene observes, did not facilitate its comparison with
classical mathematics.  Kleene and Vesley's book ``The Foundations of
Intuitionistic Mathematics'' \cite{FIM}, which was published in 1965, was
the result of many years' research.  It presented a formal system $\mathcal{FIM}$
for a part of intuitionistic mathematics including Brouwer's set theory restricted
to definable properties of numbers and sequences, sometimes referred to as
intuitionistic analysis.\footnote{Kleene himself referred to this part as ``two-sorted
intuitionistic number theory.''}  The language he uses is identical with that of
a formal system for classical analysis.  However, in contrast to intuitionistic
propositional and predicate logic and Heyting arithmetic, intuitionistic analysis
\textit{is not} a subsystem of the classical theory, but a divergent theory of
analysis, as would be expected from the discussion above.  However,
$\mathcal{FIM}$ contains a subsystem, the \textit{basic system} \,$\mathcal{B}$,
consisting of the common part of classical and intuitionistic analysis.
When X10$^c$ (the principle of double negation) is added to \,$\mathcal{B}$ the
result is a corresponding system for classical analysis, while the intuitionistic
system $\mathcal{FIM}$ comes from $\mathcal{B}$ by adding one axiom expressing
the classically false continuity principle.

\vskip 0.1cm

The language of $\mathcal{FIM}$ is a two-sorted first order language with equality,
suitable for a formal theory of natural numbers and choice sequences (functions from
$\omega $ to $\omega $); we shall see that it can also express continuous functionals.
There are two types of individual variables, arithmetical variables a,b,c,\ldots ,
x,y,a, \ldots and function variables $\alpha,\beta,\gamma,\ldots $, and finitely
many function constants $\mathrm{f_{0},\ldots,f_{24} ,}$ expressing zero, successor,
addition and multiplication, as well as other primitive recursive functions such as
predecessor, remainder and quotient; each $\mathrm{f_{i}}$ expresses a function
of $k_i$ number and $l_i$ function arguments, respectively.  This choice is not
unique, it was intended by Kleene to simplify the development of the
theory.\footnote{The complete catalog of function constants can be found in \cite{FIM}.
In the spirit of Brouwer it is open, allowing additions as required by the development
of the theory; thus, in \cite{RM} two more function constants are proposed.}
In addition, the symbols of the language include Church's $\lambda$.

\vskip 0.1cm

In addition to the \textit{terms}, which are the formal expressions of
$\mathcal{FIM}$ for the natural numbers, in this system there are formal
expressions, called \textit{functors}, for (total) functions from $\omega $
to $\omega$.  The two notions are defined by simultaneous induction:
\begin{enumerate}
\item[(i)] The number variables are \textit{terms}.\\[-10pt]
 \item[(ii)] The function variables are \textit{functors}.\\[-10pt]
  \item[(iii)] For every \;$i=0,\ldots,24$, if \;$k_i=1$ \;and \;$l_i=0$, then
\;$\text{f}_i$\; is a \textit{functor}.\\[-10pt]
 \item[(iv)] For every \;$i=0,\ldots,24$, if \;
$\text{t}_1,\ldots,\text{t}_{k_i}$ \;are \textit{terms} and
$\text{u}_1,\ldots,\text{u}_{l_i}$ \;are \textit{functors}, then
\; $\text{f}_i(\text{t}_1,\ldots,\,\text{t}_{k_i},\text{u}_1,\ldots,
\text{u}_{l_i})$ \; is a \textit{term}.\\[-10pt]
  \item[(v)] If \;u \;is a \textit{functor} and \;t is a \;\textit{term}, then
\;(u)(t)\; is a \textit{term}.\\[-10pt]
 \item[(vi)] If \;x \;is a number variable and \;t is a \;\textit{term}, then
\;$\lambda\text{x\,t}$ \;is a  \textit{functor}.\\[-10pt]
 \item[(vii)] An expression is a \textit{term} or \textit{functor} if and only
if that follows from (i)--\;(vi).
\end{enumerate}
The atomic or prime formulas are expressions of the form s = t where
s and t are terms.  Equality between functors is not prime, but is defined by
$$\mathrm{u=v\;\equiv \;\forall x \;(u(x)=v(x))},$$ where \;x\; does not appear free
in \;u, v.  The \textit{formulas} are formed as usual, with the difference that
here the quantifiers may apply to function variables, so that
$\forall \text{áA,}\;\;\exists\text{áA}\;\;$ are formulas.

\vskip 0.1cm

The free and bound occurrences of variables are defined as usual, except that
Church's $\lambda$ symbol also binds number variables.

\vskip 0.1cm

The axiom schemas and rules of inference of the system include first of all
X1 - X12 and R1 - R3 (for the language of $\mathcal{FIM}$) as well as the
following, where u is a functor free for $\alpha$ in A($\alpha$) and $\alpha$
does not occur free in C:
\begin{enumerate}
\item[X13.]{$\forall$á A(á) $\rightarrow $ A(u).}
\item[X14.]{A(u) $\rightarrow \exists á A(á)$.}
\item[R4.]{From C $\rightarrow$ A(á), conclude C $\rightarrow \forall$á A(á).}
\item[R5.]{From A(á) $\rightarrow$ C, conclude $\exists$á A(á) $\rightarrow$ C.}
\end{enumerate}
The equality axioms ×Å1 - ×Å3 and the arithmetical axioms ×Ind and ×Í1 - ×Í7
are also included, where of course x, y, z are number variables. The axioms
concerning functors are first of all the definitions of the function constants
$\mathrm{f_{0},\ldots,f_{24},}$ which are formulas expressing the corresponding
explicit or recursive definitions, as well as:
 \begin{enumerate}
 \item[XF1.]{$\mathrm{(\lambda x\, r(x))(t)=r(t).}$}
 \item[XF2.]{$\mathrm{a=b\rightarrow \text{á}(a)=
  \text{á}(b).}$}
  \item[XF3.]{$\mathrm{\forall x \,\exists \text{á}
  \,A(x,\text{á})\rightarrow \exists\text{á}\,\forall x \,A(x,\lambda
  y\,
  \text{á}(2^{x}\cdot 3^{y})),}$}
  \end{enumerate}
where \;r(x), t \;are terms, \;x\; is a number variable, t is free for x in \;r(x),
and Á(x,\,á) is a formula in which x is free for á.

XF1 is the $\beta$-reduction axiom of the ë-calculus; XF3 is a principle of
countable choice which is acceptable according to the intuitionistic
interpretation of the quantifiers.

\vskip 0.1cm

For the description and treatment of choice sequences in $\mathcal{FIM}$ we use
\textit{sequence numbers}, which code initial segments of sequences of natural
numbers as follows:\;\; the first \;$x$ \;values of a sequence \;$\alpha$ \;are
given by $\overline{\alpha}(x)$, where
\[\overline{\alpha}(0)=1\;\;\tx{and }\;\;
  \overline{\alpha}(x+1)=p_0^{\alpha (0)+1}\cdot\ldots\cdot
  p_{x}^{\alpha(x)+1},\]
where $p_0,p_1,\ldots$  are the prime numbers in their natural order.
This coding, as well as the formal predicate Seq(a) expressing that \;a
\;codes a finite sequence, can be defined in the formal system.  So can
the operation  $\mr{a*b=a\cdot\prod_{i<lh(b)}p_{lh(a)+i}^{(b)_i}}$\footnote{Here
lh(a) expresses the length of the sequence coded by a, and $\mr{(a)_i}$ its
i$^\mr{th}$ projection.} of concatenation.
With these tools it becomes possible to formulate in $\mathcal{FIM}$ the
characteristic principles of Brouwer's mathematics concerning choice sequences.

The mathematical principle which Brouwer tried to convey by the bar theorem
is introduced in $\mathcal{FIM}$ by the Axiom of Backwards Induction
(Bar Induction):

\vskip 0.2cm

 $\begin{array}{ll}
 \mathcal{BI}:\;\;\;\; &\mathrm{\forall a \;[\,Seq(a)\rightarrow R(a)\vee \neg
 R(a)\,]\; ~\&~ \;\forall } \text{á} \,\exists\mathrm{x\,
 R}(\overline{\text{á}}(\mathrm{x))\,~\&~ }\;\s
 &\forall\mathrm{ a\; [\,Seq(a)~\&~ R(a)\rightarrow A(a)\,]\,~\&~}\s
 &\forall\mathrm{ a\; [\,Seq(a)~\&~ \forall s\, A(a*2^{s+1})\rightarrow
 A(a)\,]}\s
  &\rightarrow \mathrm{A(1).}
  \end{array}$

\vskip 0.2cm

About Kleene's choice  of this schema we can say the following:\\
 The properties of choice sequences which are interesting from the
  intuitionistic viewpoint are those which are verifiable on the basis of
  some initial part of each choice sequence, that is they are of the form
  \;$\exists x R(\overline{\alpha}(x))$, \;where \;$R(a)$ \;is an arithmetical
  predicate, decidable at least for sequence numbers.  $\mathcal{BI}$ expresses
  the following induction principle for spreads (trees) on which there is a
  (decidable) bar:  if
(i) \;$A(a)$ \, is a property which is possessed by sequence numbers having the
  property \;$R(a)$ \;and
(ii) \;$A(a)$ \,is transmitted\footnote{In the direction from the bar
  \textit{toward} the root of the tree, hence the name \textit{backwards}
  induction.} recursively, in the sense that for every sequence
  number \;$a$, \;if  \;$Á(a*2^{s+1})$\; holds for every \;$s$, then also
  \,$A(a)$ holds,
then we conclude that \,$A(1)$ holds, where 1 is the code of the empty sequence.\\
$R(a)$ \;plays the role of the bar; this form of the axiom only takes care of the
case of the universal spread, but the notion of spread is definable in $\mathcal{FIM}$
and the general form is provable from this.

The important Fan Theorem, formally stated, is proved in  $\mathcal{FIM}$ using
$\mathcal{BI}$.

$\mathcal{BI}$ is provable classically (and without the restriction of
decidability of \;$R(a)$, which is meaningless in this case); together with
the preceding axioms and rules, it completes the basic system $\mathcal{B}$.

\vskip 0.1cm

The next and final axiom expresses the classically false continuity principle.
For the formulation of this principle in $\mathcal{FIM}$ Kleene points out that
Brouwer himself, speaking about the assignment of a value $b$ to a
function $\alpha$, used the expression ``\textit{the algorithm of the
correlation law}''.  Taking this expresssion seriously and arguing as follows:

        first, the algorithm for computing the value $b$ must
        (i) decide, for each initial segment $\alpha(0),\ldots ,\alpha(y-1)$
        of a choice sequence $\alpha$, whether these values suffice for the
        computation of $b$;
        (ii) if the answer is ``yes,'' compute $b$, and

        second, this information must be given by a function $\tau$,
        which acts on sequence numbers $\overline{\alpha}(y)$ and
        produces 0 as long as the algorithm does not answer ``yes,''
        while if $y$ is the least number for which the answer is ``yes''
        then $\tau (\overline{\alpha}(y))=b+1$ (and we may assume without
        loss of generality that this is the only $y$ for which
        $\tau (\overline{\alpha}(y)) > 0$),

\vskip 0.2cm

\noindent Kleene proposed the following formal statement of the continuity
principle (Brouwer's principle):

\vskip 0.2cm

        $\begin{array}{ll}
         \mathcal{BP}:\;\; &
         \forall\text{á}\,\exists\text{â}\,\mathrm{A}(\text{á,â})
         \rightarrow\exists\text{ô}\forall\text{á}\,[\forall\mathrm{t\,
         \exists !y}\,\text{ô}(2^{\mathrm{t+1}}*\overline{\text{á}}(
         \mathrm{y}))>0\,~\&~\s
         & \forall\text{â}\,[\forall\mathrm{t\,\exists y}\,\text{ô}(2^{
         \mathrm{t+1}}*\overline{\text{á}}(\mathrm{y}))=\text{â}(
         \mathrm{t)+1\rightarrow A}(\text{á,â})\,]\,].
         \end{array}$

\vskip 0.2cm

This principle is false for classical mathematics.  For example, in
\cite{FIM} it is used to prove
$$\vdash\neg\forall\text{á}(\forall\text{x\;á(x)}=0 \vee\neg\forall\text{x\;á(x)}=0)$$
(from which it follows that equality between choice sequences must be undecidable),
and the negation of the universal closure of the least number principle.

In addition, the third chapter of \cite{FIM} contains R. E. Vesley's presentation
of the intuitionistic theory of the continuum in the context of the formal system
$\mathcal{FIM}$.  Real numbers are defined by formalizating the notion of real
number generator and the uniform continuity theorem is proved.  In the fourth
and last chapter, Kleene treats formally the question of the order of the reals.
Thus an important part of Brouwer's mathematics is faithfully represented in
this formal system.

\subsection{Function realizability interpretation of the system $\mathcal{FIM}$}

In the same book, by adapting the basic idea of number realizability, Kleene
gives an interpretation of intuitionistic analysis which guarantees the
consistency of the formal system.  Here we have propositions of the form
\;$\exists\beta\,A(\beta)$, \;where \;$\beta$ is a function variable; to interpret
and justify such a proposition requires a function \;$\beta$ \; and a verification
that the interpretation of \;$A(\beta)$ is correct.  Now the suitable candidates
for realizing objects are not natural numbers, but \textit{number-theoretic
functions of one (number) variable}.  The basic mechanism Kleene used, in
analogy with the coding of recursive functions by natural numbers, is the coding
of (continuous) functionals $F: \,\omega^{\omega}\rightarrow \omega^{\omega}$
by means of functions $\tau:\,\omega\rightarrow\omega$.  Thus he defined the
partial recursive function $\{\tau\}[\alpha]$ of two function variables $\tau$
and $\alpha$, which computes the functional coded by $\tau$ at the argument $\alpha$.

\vskip 0.1cm

The definition of function realizability is recursive and is motivated
by the previous discussion.\footnote{It also uses the definitions
$(\epsilon )_i=\textit{ë} t\;(\epsilon (t))_i$ for i = 0,1,
and $\{\epsilon\}[x]=\{\epsilon\}[ \textit{ë} t \;x]$.}  The concept to be
defined is
  $$\epsilon \;r_Ø\;\mathrm{E}:\;
    \; \epsilon \;\text{realizes}\;\mathrm{E}\;
    \text{with respect to} \;Ø,$$
    where \;$\epsilon :\omega\rightarrow\omega$,\;E \;is a formula, \;Ø\; is a
    list of variables containing all those occurring free in \;E,\; and \;$Ø$\;
    is an assignment of numbers and (one-place) number theoretic functions to
    the variables of \;Ø.\footnote{We read $\downarrow$ as ``is defined''.}
\begin{enumerate}
\item[1.] $\epsilon\;r_Ø\;\mathrm{P}$,\;\;where \;P\; is prime\;
    $\Longleftrightarrow_{\texttt{df}}$\;P is true for $Ø$,\\
    that is, P holds for the values $Ø$ assigns to the variables of Ø.
    \hfill $\mathcal{R}_P$\\[-8pt]
\item[2.] $\epsilon \;r_Ø\;\mathrm{A~\&~
B}\;\;\;\Longleftrightarrow_{\!\texttt{df}}\;(\epsilon
)_0\;r_Ø\;\mathrm{A}\;\;\text{and}\;\;(\epsilon
)_1\;r_Ø\;\mathrm{B}.$ \hfill $\&\:\mathcal{R}$\\[-8pt]
\item[3.] $\epsilon \;r_Ø\;\mathrm{A\vee B\;\;\;\Longleftrightarrow_{\!
\texttt{df}}}\;\; \text{if}\;(\epsilon
(0))_0=0\;\text{then}\;(\epsilon )_1\;r_Ø\;
\mathrm{A}\;\text{and}\\
\hspace*{105pt}
 \text{if}\;(\epsilon (0))_0\neq 0\;\text{then}\;(\epsilon
)_1\;r_Ø\; \mathrm{B}.\hfill \vee\mathcal{R}$\\[-8pt]
\item[4.] $\epsilon \;r_Ø\;\mathrm{A\rightarrow B\,\:
\Longleftrightarrow_{\! \texttt{df}}}\;\; \text{for every}\;\alpha
\;(:\omega\rightarrow\omega),\;
\text{if}\;\, \alpha\; r_Ø\;\mathrm{A}\;\,\text{then}\;\\
\hspace*{104pt}
 \{\epsilon\}[\alpha
]\!\downarrow\;\text{and}\;\{\epsilon\}[\alpha ]\;
r_Ø\;\mathrm{B}.\hfill\rightarrow\!\!\mathcal{R}$\\[-8pt]
\item[5.] $\epsilon\;r_Ø\;\neg\mathrm{A}\quad\quad
\Longleftrightarrow_{\! \texttt{df}}\; \text{for every}\;\alpha,\;
\text{it is not the case that }\;\;\alpha\;r_Ø\;\mathrm{A}$\\
\hspace*{101pt} (equivalently, if and only if \;$\epsilon\;r_Ø\;(\mathrm{A}
\rightarrow 1=0)$,\\ \hspace*{101pt}
 because of 4 and 1: there is no\;$\alpha$\; such that
 \hfill$\neg\mathcal{R}$\\
 \hspace*{101pt}
 $\alpha\;r_Ø\;1\!=\!0$, where 1\,=\,0 is a false\\
 \hspace*{101pt} prime formula).\\[-8pt]
 \item[6.] $\epsilon\;r_Ø\;\forall\mathrm{x\,A}\;\quad\,
 \Longleftrightarrow_{\! \texttt{df}}\; \text{for every}\;x\;(\in
 \omega ),\;\{\epsilon\}[x]\!\downarrow \;\;\text{and} \\
 \hspace*{103pt}\{\epsilon\}[x]\;r_{Ø,\,x}\;\mathrm{A},\;
 \text{where}\,\;x\;\,\text{is the value of \; x}.\hfill\forall
 N\mathcal{R}$\\[-8pt]
 \item[7.] $\epsilon\;r_Ø\;\exists\mathrm{x\,A}\;\quad\,
 \Longleftrightarrow_{\! \texttt{df}}\;(\epsilon )_1\;
 r_{Ø,\,(\epsilon (0))_0}\;\mathrm{A},\;\;
 \text{where}\;\; (\epsilon (0))_0 \\
 \hspace*{105pt}
 \text{is the value of \;x}.\hfill
 \exists N\mathcal{R}$\\[-8pt]
 \item[8.] $\epsilon\;r_Ø\;\forall\text{á}\,\mathrm{A}\;\quad\,
 \Longleftrightarrow_{\! \texttt{df}}\; \text{for every}\;\alpha,\;
 \{\epsilon\}[\alpha]\!\downarrow\;\text{and}\;\{\epsilon\}[\alpha]
 \;r_{Ø,\,\alpha}\;\mathrm{A},\\
 \hspace*{105pt}
 \text{where}\;\;\alpha\;\; \tx{is the value of \;á}.
 \hfill\forall F\mathcal{R}$\\[-8pt]
 \item[9.] $\epsilon\;r_Ø\;\exists\text{á}\,\mathrm{A}\;\quad\,
 \Longleftrightarrow_{\! \texttt{df}}\;\{(\epsilon )_0\}
 \!\downarrow\;\;\text{and}\;\;\;(\epsilon )_1\;r_{Ø,\{(\epsilon)_0\}}
 \;\mathrm{A},\\
 \hspace*{105pt}
 \text{where}\;\; \{(\epsilon )_0\}\;\;\text{is the value of \;á}.
 \hfill\exists F\mathcal{R}$\\[-8pt]
\end{enumerate}

A \textit{closed} formula \;E \;is\; \textit{realizable}, if it is realized
by some \textit{general recursive function} \;$\epsilon:\omega\rightarrow\omega$.
An \textit{open} formula is \textit{realizable}, if its universal closure is.

\vskip 0.1cm

The appropriate soundness theorem holds:\\
 \textit{Theorem.} (Kleene).  Let \;Ã\; be a (finite) list of formulas and
 \;Å\; a formula.  Then, if \;Ã $\vdash_{\mathcal{FIM}}$ Å \;\;and the formulas
 of  \;Ã \;are realizable, it follows that \;E \;is realizable.

\vskip 0.1cm

\textit{Corollary.} The formal system $\mathcal{FIM}$ is consistent.

\vskip 0.1cm

Completing the presentation of function realizability, Kleene made the
conjecture that, with the formalization of this notion and of the
(model theoretic) proof of the consistency of $\mathcal{FIM}$, a
metamathematical proof of the relative consistency of $\mathcal{FIM}$ with
respect to the basic system $\mathcal{B}$ would be obtained.  In 1969,
in his monograph ``Formalized recursive functionals and formalized
realizability'' \cite{RN}, he presented a detailed formalization of
the theory of recursive functions of type 2 and the corresponding
formal notion of function realizability\footnote{As in the arithmetical
case, what is involved is a translation, by which to each formula E there
corresponds a formula $\varepsilon$ r E.} as well as a variation which
made it possible to prove his conjecture, and also to prove the disjunction
and existence properties and a form of Church's thesis, for $\mathcal{FIM}$.

\vskip 0.1cm

\textit{Corollary.}  (i) If $\vdash_\mathcal{FIM}\mr{A\vee B}\;$ where
      $\mr{A\vee B}$ is a closed formula, then one of the following holds:
      $\; \vdash_\mathcal{FIM}\mr{A}$\;\;or\;\;\;$\vdash_\mathcal{FIM}\mr{B}$.

      (ii)  If \;$\vdash_\mathcal{FIM}\;\mr{\exists x\,A(x)}$ where
      \;$\mr{\exists x\,A(x)}$\; is a closed formula, then, for some natural
      number $x$, \;$\vdash_\mathcal{FIM}\;\mr{A(\mathbf{x})},$ where
      $\mathbf{x}$ is the numeral for $x$.

      (iii)  If $\vdash_\mathcal{FIM}\mr{\exists\text{á}\,A(\text{á})}$ where
      $\mr{\exists \text{á}\,A(\text{á})}$ is a closed formula, then
      $\vdash_\mathcal{FIM}\mr{\exists\text{á}_{GR(\text{á})} A(\text{á})}$,
      where GR(á) is a formula (with only á free) expressing that the function
      represented by á is general recursive.

\vskip 0.1cm

In 1973 Troelstra characterized function realizability using the schema GC$_1$,
a generalization of the continuity principle $\mathcal{BP}$.\footnote{GC$_1$ is
the schema $\mathrm{\forall \tx{á} (A(\tx{á})\rightarrow\exists\tx{â} B(\tx{á}
,\tx{â}))\rightarrow\exists\tx{ó}\forall\tx{á}
(A(\tx{á})\rightarrow\exists\tx{â} (\{\sigma\}[\tx{á}]=\tx{â}~\&~ B(\tx{á},\tx{â})))}$,
where A is an almost negative formula.}

\subsection{Relativized and modified realizabilities}
Kleene relativized his original function-realizability interpretation in various ways.
If $\Phi$ is a class of number-theoretic functions closed under ``recursive in'' then
for $~\Theta, \Psi \subseteq \Phi~$ and $\varepsilon \in \Phi$ the notion
``$\varepsilon \,$ $\Phi$/{\em realizes}-$\Psi \; \mathrm{E}$'' is defined like
``$\varepsilon \,$ realizes-$\Psi \; \mathrm{E}$'' except that in Clauses 4, 5 and 8 the
$\alpha$ is restricted to $\Phi$.

Then $\mathrm{E}$ is $~\Phi$/{\em realizable}/$\Theta~$ if for some $\varepsilon$
recursive in $\Theta$:  $~\varepsilon \,$  $\Phi$/realizes $\, \mathrm{\forall E}$.
The soundness theorem extends to the relativized notions.  Using his example of a
binary fan all of whose recursive branches (but not all of whose branches) are
finite, and taking $\Phi = \Theta$ to be the class of all general recursive
functions, Kleene showed that Brouwer could not have proved his ``Bar Theorem''
without using bar induction in the proof.

\vskip 0.1cm

{\em Theorem}.  (Kleene)  The axiom schema of bar induction is independent of the
other axioms of intuitionistic analysis.  If the ``fan theorem'' replaces the bar
induction schema, then it too is independent of the other axioms.\footnote{\cite{FIM}
Corollary 9.9.  The circularity in Brouwer's ``proof'' of his Bar Theorem is analyzed
in Chapter 6.}

\vskip 0.1cm

Next he took $\Phi = \Theta$ to be the class $\Xi$ of all arithmetical functions,
which is (classically) closed under general recursiveness and the jump operation,
and proved (classically) that the fan theorem and all axioms of the intuitionistic
system {\em except} bar induction are $~\Xi$/{\em realizable}/$\Xi~$.  (Later,
Howard and Kreisel \cite{HK1966} proved that bar induction actually changes
intuitionistic arithmetic; Troelstra \cite{Tr1974} proved that the fan theorem
does not.)

\vskip 0.1cm

Kreisel \cite{Kr1958} first suggested a different kind of modification of
realizability (later adapted by Kleene) in order to prove Markov's Principle
independent of the intuitionistic axioms.  Kreisel and Kleene used explicit
types, but the same effect is obtained using implicit types via a notion of
{\em agreement}.  Van Oosten \cite{VO2000} explains the main idea:
``Each formula gets {\em two} sets of realizers, the {\em actual} realizers
being a subset of the {\em potential} ones.''

\vskip 0.1cm

For example, $\varepsilon$ {\em agrees with} $\mathrm{A \rightarrow B}$ if,
whenever $\sigma$ {\em agrees with} $\mathrm{A}$, $~\{\varepsilon\}[\sigma]$ is
defined and {\em agrees with} $\mathrm{B}$.  If {\sf F} is a collection of functions
closed under ``general recursive in'' and the free variables of $\mathrm{A \rightarrow B}$
are interpreted by numbers and functions from {\sf F}, then
$\varepsilon$ $^{\sf F}${\em realizes} $\mathrm{A \rightarrow B}$ under this
interpretation of the variables if $\varepsilon \in {\sf F}$ and $\varepsilon$ agrees
with $\mathrm{A \rightarrow B}$ and for every $\sigma \in {\sf F}$:
if $\sigma$ $^{\sf F}${\em realizes} $\mathrm{A}$ under the interpretation then
$\{\varepsilon\}[\sigma]$ $^{\sf F}${\em realizes} $\mathrm{B}$ under the interpretation.

\vskip 0.1cm

The proof that bar induction is $^{\sf F}$realizable is a little more complicated
than for realizability.  The critical observation is that if $\sigma~$
$^{\sf F}$realizes $\mathrm{\forall \alpha \exists ! x R(\overline{\alpha}(x))}$, then
$\{\sigma\}[\varphi]$ is defined for {\em every} $\varphi$ by agreement, and
$(\{\sigma\}[\varphi](0))_0 = x$ depends only on a finite initial segment of
$\varphi$.  But every neighborhood of $\varphi$ contains elements of {\sf F},
and the assumption on $\sigma$ implies (by a roundabout argument) that if
$\overline{\psi}(x) = \overline{\varphi}(x)$ and $\psi \in {\sf F}$ then
$(\{\sigma\}[\psi](0))_0 = x$ also.  This provides the basis for informal
bar induction.

\vskip 0.1cm

For Kleene's $_{\sf s}$realizability, {\sf F} is the class of all functions.
To see that Markov's Principle is not $_{\sf s}$realizable, assume
$\sigma$ $_{\sf s}$realizes $\mathrm{\forall \tx{á} (\neg \forall x \neg (\tx{á}(x) = 0)
\rightarrow \exists x (\tx{á}(x) = 0))}$.  Then $\tau = \{\sigma\}[\lambda t. 1]$
$_{\sf S}$realizes
$\mathrm{\neg \forall x \neg (\tx{á}(x) = 0) \rightarrow \exists x (\tx{á}(x) = 0)}$ when
$\tx{á}$ is interpreted by $\lambda t. 1$, and there is a recursive $\theta$
which agrees with the hypothesis, so $(\{\tau\}[\theta](0))_0 = x$ is determined by
$\overline{\lambda t. 1}(z)$ for some $z$.  If $\varphi(y) = 1$ for $y \leq x + z$
but $\varphi(x + z + 1) = 0$, then $\theta$ $_{\sf s}$realizes the hypothesis
$\mathrm{\neg \forall x \neg (\tx{á}(x) = 0)}$ when $\tx{á}$ is interpreted by $\varphi$,
but $\{\tau\}[\theta]$ does not $_{\sf s}$realize the conclusion.

\vskip 0.1cm

The relativized version was developed in \cite{JRM1971} (for {\sf F} the class of all
recursive functions) in order to prove $\mathrm{\forall \tx{á} \neg \neg GR(\tx{á})}$
consistent with the intuitionistic theory and with Vesley's Schema \cite{Ve1970},
which entails Brouwer's creating subject counterexamples.  For other applications,
first observe that ordinary function-realizability provides a classical proof that
the intuitionistic theory {\bf FIM} is consistent with {\bf PA}, since classically
every {\em arithmetical} instance of the law of excluded middle is realized by
some function.

\vskip 0.1cm

If $\mathrm{A(x)}$ is any arithmetical predicate with only $x$ free, let
{\sf F}[$\mathrm{A(x)}$] be the collection of all functions classically recursive
in the intended interpretation of $\mathrm{A(x)}$.  Then both
$\mathrm{\forall \tx{á}  \neg \neg \exists e \exists \tx{â}
[\forall x (\tx{â} (x) = 0 \leftrightarrow A(x))~\&~
\forall x \exists y [T_1(e,x,\overline{\tx{â}}(y))~\&~
U(\overline{\tx{â}}(y)) = \tx{á}(x)]]} \,$ and $\mathrm{\forall x (A(x) \vee \neg A(x))}$
are $^{{\sf F}[\mathrm{A(x)}]}$realizable and hence
consistent with {\bf FIM}.  Thus it is possible to interpret the constructively
undetermined part of the intuitionistic continuum as classically recursive in
the constructively determined part, which could include the characteristic
functions of arbitrarily complicated arithmetical relations.  Of course,
it is impossible to assign g\"odel numbers (or even relative g\"odel numbers)
continuously, so the $\mathrm{\neg \neg \exists e}$ cannot consistently be replaced by
$\mathrm{\exists e}$.

\subsection{A glimpse at today and tomorrow}

His philosophical ideas were Brouwer's motivation for the intuitionistic
reconstruction of mathematics.  However, in an irony of history, the
interest in intuitionistic thought (except among those holding similar
philosophical ideas) springs today mostly from logic and computer science.

\vskip 0.1cm

The B-H-K interpretation found a precise implementation in Kleene's
notion of realizability, which relates it to the theory of recursive
functions; but also in connection with the proof systems of Gentzen, it
led to the Curry-Howard isomorphism, which correlates ë-terms and finally
programs with proofs.  The intuitionistic theory of types of Per Martin-L\"{o}f
belongs in the same framework.  The polymorphic ë-calculus or system F of
Girard and the logic of constructions of Coquand and Huet (on which was based
the functional language proof checker Coq, used recently (2004) to verify
the correctness of the solution to the 4 color problem), are applications
of the preceding.  A detailed description of the development and further
influence of these ideas appears in the article ``From constructivism to
computer science'' by A. S. Troelstra, \cite{Tr1999}.

\vskip 0.1cm

Another source of interest is category theory, where it was discovered that
important spaces such as topoi are models of intuitionistic logic.  In the
article \cite{VO2000} by J. van Oosten one can find a sketch of this line of
development of the theory.

\vskip 0.1cm

Various logical and mathematical systems and semantics were developed
in connection with intuitionistic logic.  As examples we mention
Kripke semantics and constructive versions of Zermelo - Fraenkel
set theory, like those of J. Myhill and P. Aczel; also the constructive
development of a large part of analysis by E. Bishop and his school, and
the recursive constructive mathematics of the Russian school associated with
A. A. Markov.

\vskip 0.1cm

In each case, intuitionistic mathematics constitutes a rich source of
mathematical ideas.  These can either be developed within classical mathematics
to produce proper subtheories of particular classical theories, or be
considered as extending classical mathematics,\footnote{As specific cases of
the mathematical implementation of this viewpoint we mention
\cite{JRM1996, JRM2003}} offering different definitions and refinements
of logical and mathematical concepts, with all that this implies.

\end{document}